\documentclass{svmult}

\usepackage{mathptmx}       
\usepackage{helvet}         
\usepackage{courier}        
\usepackage{type1cm}

\usepackage{amsmath,amssymb,bm} 
\usepackage{epsfig,color,graphicx}
\usepackage{subfigure}
\usepackage{enumerate}

\usepackage[group-separator={,}]{siunitx}

\usepackage[colorlinks=true,linkcolor=black,citecolor=black,urlcolor=black,bookmarks,breaklinks=true]{hyperref}

\usepackage{empheq}
\usepackage{mathtools}
\mathtoolsset{showonlyrefs=true}

\newcommand{\1}{{\mathbf{1}}} 

\def\S{{\mathbb{S}}}

\renewcommand{\E}{{\mathbf{E}}}
\renewcommand{\P}{{\mathbf{P}}}

\def\R{{\mathbb{R}}}
\def\N{\mathbb{N}}
\def\cK{{\mathcal K}}
\def\cR{{\mathcal R}}
\def\cA{{\mathcal A}}
\def\cH{{\mathcal H}}
\def\overpsi{\overline{\Psi}}

\DeclareMathOperator{\cov}{\mathbf{cov}} 
\DeclareMathOperator{\var}{\mathbf{v}\mathbf{a}\mathbf{r}} 

\newcommand{\BQ}{\mathbf{Q}}

\begin{document}

\title*{Second order analysis of geometric functionals of Boolean models}
\author{Daniel Hug, Michael A. Klatt, G\"unter Last, Matthias Schulte}
\institute{Daniel Hug \at Karlsruhe Institute of Technology, Institute of Stochastics, D-76128 Karlsruhe, 
\email{daniel.hug@kit.edu}
\and Michael A. Klatt \at Karlsruhe Institute of Technology, Institute of Stochastics, D-76128 Karlsruhe, 
\email{michael.klatt@kit.edu}
\and G\"unter Last \at Karlsruhe Institute of Technology, Institute of Stochastics, D-76128 Karlsruhe, 
\email{guenter.last@kit.edu}
\and Matthias Schulte \at Karlsruhe Institute of Technology, Institute of Stochastics, D-76128 Karlsruhe, 
\email{matthias.schulte@kit.edu}}
\maketitle

\abstract{
This paper presents asymptotic covariance formulae and central limit
theorems for geometric functionals, including volume, surface area, and
all Minkowski functionals and translation invariant Minkowski tensors as
prominent examples, of stationary Boolean models. Special focus is put
on the anisotropic case.
In the (anisotropic) example of aligned rectangles, we provide explicit
analytic formulae and compare them with simulation results.
{We discuss which information about the grain distribution second
moments add to the mean values.}
}

\section{Introduction}\label{sec:Introduction}

In this article we study a large class of functionals of the Boolean
model, a fundamental benchmark model  of stochastic
geometry~\cite{CSKM13, SchneiderWeil, Molchanov1996} and continuum
percolation~\cite{Hall1988, MeesterRoy96}. It has many applications in
material science~\cite{Torquato.2002}, physics~\cite{Arns.2003,
ScholzEtAl2015}, and astronomy~\cite{Buchert.1994, Kerscher.2001}, as
well as, for the measurement of biometrical data~\cite{Mecke.2005} or
the estimation of percolation thresholds~\cite{Mecke.1991, Mecke.2002b}.
Intuitively speaking, a Boolean model is a collection of overlapping
random grains, scattered in space in a purely random manner. This random
object is defined as follows.  Let $X=\{X_1,X_2,\hdots\}$ be a
stationary Poisson process of intensity $\gamma$ in $\R^n$, that is, a
countable collection   of random points in $\R^n$ such that the numbers
of points in disjoint sets are
independent and the number of points in each set follows a Poisson
distribution whose parameter is $\gamma$ times the Lebesgue measure of
the set. Let $(Z_i)_{i\in\N}$ be a
sequence of independent and identically distributed random convex bodies
(nonempty compact convex subsets of $\R^n$), independent of $X$. The
Boolean model $Z$ is the
random closed set defined by
$$
Z:=\bigcup_{i\in\N} (Z_i+X_i),
$$
where $Z_i+X_i:=\{z+X_i:z\in Z_i\}$. An example is the spherical Boolean
model, where the $Z_i$ are balls with random radii centred at the
origin and $Z_i+X_i$ is the corresponding ball centred at $X_i$.

In this paper we study  geometric functionals of the Boolean model $Z$.
Prominent examples of such functionals are the intrinsic volumes
(Minkowski functionals) and Minkowski tensors, which are efficient shape
descriptors that have been successfully applied to a variety of physical
systems~\cite{SchroederTurketal2011AdvMater}.
(In~\cite{SchroederTurketal2013NJP}, the different approaches and
notations in the physics and mathematics literature are compared.).  We
are interested in the second order properties of the random variables
obtained by applying geometric functionals to the restriction $Z\cap W$
of the Boolean model to a convex observation window $W\subset \R^n$. For
a stationary and isotropic Boolean model, Miles~\cite{Miles76} and
Davy~\cite{Davy76} obtained explicit formulae expressing the  mean
values of the Minkowski functionals in terms of the intensity and
geometric mean values of the typical grain (see
also~\cite{CSKM13,SchneiderWeil}). For mean value formulae for more
general functionals in Boolean models we refer
to~\cite{HoerrmannWeil2016}.
We shall discuss here formulae for
asymptotic covariances as well as multivariate central limit theorems
for an increasing observation window.  Much of the presented theory is
taken from~\cite{HLS}.  However, some results are new. In particular
this is the first paper providing explicit covariance formulae involving
the Euler characteristic of planar non-isotropic Boolean models.  Our
methods are based on the Fock space representation of Poisson
functionals from~\cite{LaPe11} and the Stein-Malliavin approach to their
normal approximation~\cite{Peccatietal2010,LPS16,LastPenrose}.  A
completely different  treatment of second moments of curvature measures
of an isotropic Boolean model with an interesting application to
morphological thermodynamics was presented in~\cite{Mecke01}. 
There, two different scenarios are considered: first, a Poisson
distributed number of grain centres in the observation window (Poisson
process), and second, a fixed number of grains (Binomial process).  In
statistical physics, these two choices are called the grand canonical
and the canonical ensemble. The second moments of geometric quantities
show a similar behaviour as thermodynamical quantities in statistical
physics~\cite{Mecke01, Mecke2000}.
For the perfectly isotropic examples of overlapping discs or spheres,
the covariances of the Minkowski functionals are also discussed in~\cite{BrMe2002} or \cite{Kerscher.2001}, respectively.

This paper is organized in the following way. After introducing Boolean
models and geometric functionals in Section~\ref{sec:Preliminaries},
Section~\ref{sec:Covariances} is devoted to the covariance structure of
geometric functionals of Boolean models. First, we present general
covariance formulae. Then, we concentrate on planar Boolean models.
Univariate and multivariate central limit theorems for geometric
functionals of Boolean models are discussed in Section~\ref{sec:CLTs}.
In Section~\ref{sec:rectangles}, we explicitly compute the covariance
formulae for a special Boolean model of aligned rectangles. 
In the final Section~\ref{sec:simulations}, we present and discuss
simulation results for Boolean models with rectangles and compare them
with our theoretical findings. The agreement is excellent.

Let us finish this introduction with an informal summary of our results
for applied scientists. We calculate for certain models of disordered
systems of overlapping grains the second moments of a quite general
class of robust shape descriptors, which include as well-known examples
volume, surface area, Euler characteristic, and, more generally, all
Minkowski functionals and tensors. Our results apply to general
anisotropic grain distributions, see
Theorems~\ref{thm:CovariancesVolumeSurfaceArea} and
\ref{thm:PositiveDefiniteness}.
The anisotropic case of aligned (planar) rectangles is discussed in
great detail; see Section~\ref{sec:rectangles} and
Fig.~\ref{fig_cov_numerical_integration}.
It is interesting to note that the asymptotic formulae for the infinite
volume system are actually exact for finite systems with periodic
boundary conditions; see Subsection~\ref{subsec:torus}.
The central limit theorem for the geometric functionals (see
Theorems~\ref{thm:MultivariateCase} and \ref{thm:UnivariateCase})
ascertains that in the limit of infinite system size the probability
distributions of the normalized geometric functionals are normal
distributions. 
If the structure of a given sample is reasonably well described by the
(joint) cumulative probability distributions of the geometric
functionals, it is possible to construct tests of certain model
hypotheses for random heterogeneous media based on the asymptotic
normality and our explicit covariance formulae. 
We discuss the behaviour of the second moments (e.g., how they differ
for various models) and probability distributions in
finite systems for specific examples: either aligned or isotropically
oriented rectangles (distributed randomly in space). Moreover, we
derive explicit formulae (see Fig.~\ref{fig_cov_numerical_integration})
and compare the results to simulations (see Figs.~\ref{fig_Boolean_cov}
and \ref{fig_Boolean_MF_epdf}).

\section{Preliminaries}\label{sec:Preliminaries}

In the introduction we have defined a Boolean model in terms of a
stationary Poisson process in $\R^n$ which is independently marked with
random convex bodies, see, e.g.,~\cite{HoerrmannWeil2016}.
In this paper we use an equivalent description based on a Poisson
process in the space $\cK^n$ of convex bodies.  For our purposes this
representation is more convenient.

We equip $\cK^n$ with its Borel $\sigma$-field $\mathcal{B}(\cK^n)$
with respect to the Hausdorff metric. We call a measure $\Theta$ on
$\cK^n$ locally finite if
$$
\Theta(\{K\in\cK^n: K\cap C\neq\emptyset\})<\infty, \quad C\in \mathcal{C}^n,
$$
where $\mathcal{C}^n$ is the space of compact subsets of $\R^n$. Let
$\mathbf{N}$ be the space of all locally finite counting measures on
$\cK^n$ and let it be equipped with the smallest $\sigma$-field
$\mathcal{N}$ such that all maps $\nu \mapsto \nu(A)$, $A\in\mathcal{B}(\cK^n)$, 
from $\mathbf{N}$ to $\N\cup\{0,\infty\}$ are
measurable. Each element $\nu\in\mathbf{N}$ has a representation
$$
\nu=\sum_{i=1}^{N} \delta_{K_i}, \quad K_1,K_2,\hdots \in\cK^n, \quad N\in\N\cup\{0,\infty\},
$$
where $\delta_K$ stands for the Dirac measure concentrated at $K\in\cK^n$. Because of this 
representation one can think of $\nu$ as a countable collection of convex bodies 
(or grains).

Throughout this paper all random objects are defined on a fixed
(sufficiently rich) probability space $(\Omega,\cA,\P)$. We call
a random element $\eta$ in $\mathbf{N}$ a Poisson process with a
locally finite intensity measure $\Theta$ if
\begin{enumerate}[(ii)]
\item[(i)] $\eta(A_1),\hdots,\eta(A_m)$ are independent for disjoint sets
  $A_1,\hdots,A_m\in\mathcal{B}(\cK^n)$,
\item[(ii)] $\eta(A)$ follows a Poisson distribution with parameter
  $\Theta(A)$ for $A\in\mathcal{B}(\cK^n)$, i.e.
$$
\P(\eta(A)=k)=\frac{\Theta(A)^k}{k!} e^{-\Theta(A)}, \quad k\in\N\cup\{0\}.
$$
\end{enumerate}
The second property explains the name. 
Since $\Theta(A)=\E \eta(A)$ for any $A\in\mathcal{B}(\cK^n)$,
$\Theta$ is called intensity measure of $\eta$. The Poisson process $\eta$ is called
stationary if it is invariant under the shifts $K\mapsto K+x:=\{y+x:y\in K\}$ for all
$x\in\R^n$. This means that the distribution of $\eta$ does not change
under simultaneous translations of its grains. The
stationarity of the Poisson process $\eta$ is equivalent to the
translation invariance of the intensity measure $\Theta$.

In the following we always assume that $\eta$
is a stationary Poisson process in $\cK^n$ with an intensity measure
$\Theta$ such that $\Theta(\cK^n)>0$. 
It follows from \cite[Theorem 4.1.1]{SchneiderWeil} that the intensity
measure $\Theta$ has the representation
$$
\Theta(\cdot) = \gamma \iint {\bf 1}\{K+x\in\cdot\} \, dx \, \mathbf{Q}(dK),
$$
where $\gamma\in(0,\infty)$ is an intensity parameter and $\mathbf{Q}$ is a probability measure on $\cK^n$ such that
\begin{equation}\label{eqn:AssumptionQ}
\int V_n(K+C) \, \mathbf{Q}(dK) < \infty, \quad C\in\mathcal{C}^n.
\end{equation}
Without loss of generality we can assume in the following that
$\mathbf{Q}$ is concentrated on the convex bodies with the centre
of the circumscribed ball as origin. A random convex body $Z_0$
distributed according to the probability measure $\mathbf{Q}$ is
called typical grain. It follows from Steiner's formula that \eqref{eqn:AssumptionQ}
is equivalent to
$$
v_i:=\E V_i(Z_0)<\infty, \quad i=0,\hdots,n,
$$
where $V_0,\hdots,V_n$ stand for the intrinsic volumes.
Later we shall require that some higher moments of the intrinsice
volumes exist. When studying covariances we have to assume that
\begin{equation}\label{eqn:AssumptionSecondMoments}
\E V_i(Z_0)^2 <\infty, \quad i=0,\hdots,n.
\end{equation}
For some results we need the stronger assumption that
\begin{equation}\label{eqn:AssumptionThirdMoments}
\E V_i(Z_0)^3 <\infty, \quad i=0,\hdots,n.
\end{equation}

The Boolean model $Z$  based on the Poisson process $\eta$
is the union of all grains of the Poisson process $\eta$, that is
$$
Z:=\bigcup_{K\in\eta} K.
$$
This is a random closed set, whose distribution is 
completely determined by the intensity $\gamma$ and the distribution of
the typical grain $Z_0$.
The stationarity of the Poisson
process $\eta$ implies the stationarity of the Boolean model $Z$,
that is, the distribution of $Z$ is invariant under
translations.
Throughout this paper we investigate the stationary Boolean model $Z$
within compact convex observation windows. For a convex body
$W\in\cK^n$ the number of convex bodies of $\eta$ that intersect $W$
is almost surely finite so that the random closed set $Z\cap W$
belongs almost surely to the convex ring $\mathcal{R}^n$, which is the
set of all unions of finitely many convex bodies and the empty set. 
Most results in this paper are for the
asymptotic regime that the observation window is increased. More
precisely, we shall assume that the inradius of the observation window goes
to infinity.

To study the behaviour of the intersection of the Boolean
model with the observation window $W$, we consider 
functionals of $Z\cap W$ with specific properties.
We say that a functional $\psi:\mathcal{R}^n\to\R$ is 
\begin{enumerate}[(iii)]
\item[(i)] additive, if $\psi(\emptyset)=0$, and
$\psi(A\cup B)=\psi(A)+\psi(B)-\psi(A\cap B)$ for all $A,B\in\mathcal{R}^n$;
\item[(ii)] locally bounded, if
$$
M(\psi):=\sup \{|\psi(K+x)|: x\in\R^n, K\in\cK^n \text{ with }K\subset [0,1]^n\}<\infty;
$$
\item[(iii)] translation invariant, if
$\psi(A+x)=\psi(A)$, for any $A\in\mathcal{R}^n$ and any  $x\in\R^n$.
\end{enumerate}
A measurable functional  $\psi: \mathcal{R}^n\to \R$ 
with all three properties is called geometric. In this case
property (ii) can be simplified using the translation invariance (iii).
Simple examples of geometric functionals are volume and surface
area. These functionals are generalized by the intrinsic volumes
$V_0,\hdots,V_n$, where $V_n$ is the volume, $V_{n-1}$ is half the
surface area (if the set is the closure of its interior) and $V_0$ is
the Euler characteristic.

More general geometric functionals are of the form
\begin{align}\label{Vgi}
V_{g,i}(A):=\Psi_i(A;g):= \int g(u) \, \Psi_i(A;du), \quad A\in\mathcal{R}^n,
\end{align}
where $\Psi_i(A;\cdot):=\Lambda_i(A;\R^n\times\cdot)$,
$i\in\{0,\hdots,n\}$, is the (additive extension of the) $i$-th
area measure of $A$ (a measure on the unit sphere $\mathbb{S}^{n-1}$
in $\R^n$), and $g:
\mathbb{S}^{n-1}\to\R$ is measurable and bounded. If $g\equiv 1$, then $V_{g,i}=V_i$. We refer to
\cite[p.~216]{Sch93} for more detail on the support measures
$\Lambda_i$.  An example for geometric functionals of the form
\eqref{Vgi} are the so-called harmonic intrinsic volumes, which are used
in \cite{DissHoerrmann} to give a representation of the intensity
$\gamma$ of non-isotropic Boolean models.

The next class of geometric functionals we consider are the components
of translation invariant Minkowski tensors
(see~\cite{Schneider2016,HugSchneider2016} for a more
detailed introduction to tensor valuations).
Let us denote by $\mathbb{T}^s$ the space of $s$-dimensional tensors in
$\R^n$. 
Let $(e_1,\ldots,e_n)$ denote the standard basis of $\R^n$. Then, for
$u\in\R^n$ and $s\in\N$, the $s$-dimensional tensor $u^s$ is given by
its coordinates
$$
(u^s)_{i_1,\hdots,i_s}=\prod_{j=1}^s u_{i_j}, \quad i_1,\hdots,i_s\in \{1,\hdots,n\}
$$
with respect to the tensor basis $e_{j_1}\otimes\cdots\otimes e_{j_s}$,
$j_1,\ldots,j_s\in\{1,\ldots,n\}$. See~\cite{HugSchneider2016} for a
description in terms of a basis of the vector space $\mathbb{T}^s$ of
symmetric tensors.

Now the Minkowski tensors $\Phi^{0,s}_m:
\mathcal{R}^n\to\mathbb{T}^s$, $s\in\N$, $m\in\{0,\hdots,n-1\}$, are
given by
$$
\Phi^{0,s}_m(A)=\frac{1}{s!}\frac{\omega_{d-m}}{\omega_{d-m+s}}
\int u^s \, \Psi_m(A;du),
$$
where $\omega_i:=i\kappa_i$ with $\kappa_i$ being the volume of
the unit ball in $\R^i$.
Each component of $\Phi^{0,s}_m$ is obviously measurable, additive and
translation invariant. For any $i_1,\hdots,i_r\in\{1,\hdots,n\}$ 
and $u\in\mathbb{S}^{n-1}$ we
have $|(u^r)_{i_1,\hdots,i_r}|\leq 1$ so that
$$
|(\Phi^{0,s}_m(K))_{i_1,\hdots,i_r}| \leq
\frac{1}{s!}\frac{\omega_{d-m}}{\omega_{d-m+s}}
\int 1 \, \Psi_m(K;du) =
\frac{1}{s!}\frac{\omega_{d-m}}{\omega_{d-m+s}} V_m(K)
$$
for $K\in\cK^n$. This shows that the components are also locally bounded.

\section{Covariance structure}\label{sec:Covariances}

We first consider general covariance formulae for geometric functionals
of Boolean models in any dimension $n$.  Then, we concentrate on planar
Boolean models and derive explicit integral formulae for the asymptotic
covariances of intrinsic volumes.

\subsection{General covariance formulae}\label{subsec:CovariancesGeneral}

In this subsection we consider the asymptotic covariance of two
geometric functionals of the Boolean model $Z$ within an observation
window $W$ as the inradius of $W$ is increased. This means that we
consider sequences of convex bodies $(W_i)_{i\in\N}$ such that
$r(W_i)\to\infty$ as $i\to\infty$, where $r(K)$ stands for the
inradius of a convex body $K\in\cK^n$. We denote this asymptotic
regime by $r(W)\to\infty$ in the sequel.

In order to present a formula for the asymptotic covariance of two
geometric functionals of a Boolean model $Z$ we have to introduce some
notation. For a geometric functional $\psi: \mathcal{R}^n \to \R$ 
the integrability assumption \eqref{eqn:AssumptionQ} implies 
for any $A\in\mathcal{R}^n$ that
$\E|\psi(Z\cap A)|<\infty$; see \cite{HLS}. 
Hence we can define $\psi^*: \mathcal{R}^n \to\R$ by
$$
\psi^*(A)=\E \psi(Z\cap A) - \psi(A), \quad A\in\mathcal{R}^n.
$$
The functional $\psi^*$ is again geometric, see~\cite[(3.11)]{HLS}.
The mapping $\psi\mapsto\psi^*$ is a key operation for the
second order analysis of the Boolean model. The following proposition
provides explicit formulae in some important examples.
To state these (and other formulae)
we need the
measure $\overpsi_{n-1}(\cdot):=\E \Psi_{n-1}(Z_0;\cdot)$.
For a bounded measurable function $g:\S^{n-1}\to \R$ we use the
notation
$$
\overpsi_{n-1}(g):=\int g(u)\,\overpsi_{n-1}(du)=\E\int g(u)\,\Psi_{n-1}(Z_0;du).
$$
The {\em volume fraction} of $Z$ is defined by $p:=\E V_n(Z\cap [0,1]^n)$
and can be expressed in the form
\begin{align}\label{p}
p=1-e^{-\gamma v_n}.
\end{align}

\begin{proposition}\label{prop1} Let $g:\S^{n-1}\to \R$ be bounded and measurable.
Then
\begin{align}\label{Vd*}
V^*_n=\; &-(1-p)V_n,\\
\label{Vd-1*}
V^*_{g,n-1}=\; &-(1-p)V_{g,n-1}+(1-p)\gamma \overpsi_{n-1}(g) V_n.
\end{align}
\end{proposition}
{\em Proof.}
Formula \eqref{Vd*} follows from an easy calculation; see \cite{HLS}.
For $j\in\{0,\ldots,n-1\}$ and $K_0\in\cK$
we obtain from Theorem 9.1.2 in \cite{SchneiderWeil} that
\begin{align*}
\E V_{g,j}(Z\cap K_0)=\sum^\infty_{k=1}\frac{(-1)^{k-1}}{k!}\int \Psi_j(K_0\cap\hdots\cap K_k;g)
\, \Theta^k(d(K_1,\hdots,K_k)).
\end{align*}
Using a result in \cite[Sections 3.2--3.4]{Hug99} or \cite[Theorem 3.1]{HoerrmannHugKlattMecke} 
(for $g\equiv 1$ see also \cite[p.~390]{SchneiderWeil}), 
we obtain that
\begin{align}\label{11}
\E V_{g,j}(Z\cap K_0)=\; &\sum^\infty_{k=1}\frac{(-1)^{k-1}\gamma^k}{k!}\\ \notag
&\times\sum^n_{\substack{m_0,\ldots,m_k=j\\ m_0+\ldots+m_k=kn+j}}
\int V^{(j)}_{m_0,\dots,m_k}(K_0,\hdots,K_k;g)
\, \BQ^k(d(K_1,\hdots,K_k)),
\end{align}
where 
$$
V^{(j)}_{m_0,\dots,m_k}(K_0,\hdots,K_k;\cdot)
:=\Lambda^{(j)}_{m_0,\dots,m_k}(K_0,\hdots,K_k;(\R^n)^{k+1}\times \cdot)
$$  
are finite Borel measures on $\mathbb{S}^{n-1}$, the mixed area measures 
of order $j$.

Consider \eqref{11} for $j=n-1$. In the summation on the right-hand side
we have $m_i=n-1$ for exactly one $i\in\{0,\hdots,k\}$ and $m_r=n$ for
$r\ne i$.
Using the decomposability
\begin{align}\label{4264}
  V^{(n-1)}_{n-1,n,\dots,n}(K_0,\hdots,K_k;g)=\Psi_{n-1}(K_0;g)V_n(K_1)\cdots V_n(K_k)
\end{align}
and the symmetry properties of the mixed area measures (see \cite{Hug99,HoerrmannHugKlattMecke}) 
we hence obtain that
\begin{align*}
&\E V_{g,n-1}(Z\cap K_0)\\
&\quad =\Psi_{n-1}(K_0;g)\sum^\infty_{k=1}\frac{(-1)^{k-1}\gamma^k}{k!}v^k_n
+V_n(K_0)\sum^\infty_{k=1}\frac{(-1)^{k-1}\gamma^k}{k!}kv^{k-1}_n \overpsi_{n-1}(g)\\
&\quad =(1-e^{-\gamma v_n})\Psi_{g,n-1}(K_0)+\gamma\overpsi_{n-1}(g)e^{-\gamma v_n}V_n(K_0).
\end{align*}
Inserting here \eqref{p} yields formula \eqref{Vd-1*}.  \hfill \hbox{}\qed

\medskip

For two geometric functionals $\psi,\phi$, we define the inner product
\begin{align}
  \varrho(\psi,\phi):=\sum_{i=1}^\infty \frac{\gamma}{i!} &\int_{\cK^n}\int_{(\cK^n)^{i-1}} 
\psi(K_1\cap \hdots \cap K_i)\label{relrho}\\
&\quad \times \phi(K_1\cap \hdots \cap K_i) \, \Theta^{i-1}(d(K_2,\hdots,K_i))
  \, \mathbf{Q}(dK_1),\nonumber
\end{align}
whenever this infinite series is well defined.
The importance of this operation for the covariance analysis of
the Boolean model is due to 
Equation \eqref{sigmaphipsi} below. In 
Proposition \ref{prop2} below and in Section \ref{subsec:CovariancesAnisotropic}
we shall see that \eqref{relrho} can be computed in some specific examples.

We need to introduce further notation.
The mean covariogram of the typical grain $Z_0$ is
$$
C_n(x)=\E V_n(Z_0\cap (Z_0+x)), \quad x\in\R^n.
$$
For a measurable and bounded function $g:\S^{n-1}\to\R$ we define
\begin{align}
C_{n-1}(x;g)=\E \int \1\{y\in Z_0^{\circ}+x\}g(u)\,\Lambda_{n-1}(Z_0;d(y,u)), \quad x\in\R^n,
\label{eq:def_C_n_minus_1}
\end{align}
where $A^{\circ}$ denotes the interior of $A$. Moreover, we use the mixed moment measures
$$
N_{n-1,n}(\cdot)=\E \iint \1\{(y,u,z)\in\cdot\}\1\{z\in Z_0\} \, \Lambda_{n-1}(Z_0;d(y,u)) \, dz
$$
and
$$
N_{n-1,n-1}(\cdot)=\E \iint \1\{(y,u,z,v)\in\cdot\} \,
\Lambda_{n-1}(Z_0;d(y,u)) \, \Lambda_{n-1}(Z_0;d(z,v)).
$$

\begin{proposition}\label{prop2} Let $g,h:\S^{n-1}\to \R$ be bounded and measurable.
Then
\begin{align}\label{VnVn}
\rho(V_n,V_n)=\; &\int \big(e^{\gamma C_n(x)}-1\big)\,dx,\\
\label{27}
\rho(V_{g,n-1},V_n)=\; &\gamma\int g(u)e^{\gamma C_n(y-z)}\,N_{n-1,n}(d(y,u,z)),\\
\label{190}
\rho(V_0,V_n)=\; &(1-p)^{-1}-1.
\end{align}
If, additionally, $\P(V_n(Z_0)>0)=1$, then
\begin{align}
\label{31}\notag
\rho(V_{g,n-1},V_{h,n-1})=\; &
\gamma^2 \int e^{\gamma C_n(y-z)}C_{n-1}(y-z;g)h(v)\,
N_{n-1,n}(d(z,v,y))\\
&\quad+\gamma \int e^{\gamma C_n(y-z)}g(u)h(v)\,N_{n-1,n-1}(d(y,u,z,v)),\\
\label{28}
\rho(V_0,V_{g,n-1})=\; &\gamma (1-p)^{-1}\overline{\Psi}_{n-1}(g).
\end{align}
\end{proposition}
{\em Proof.}
Formulae \eqref{VnVn} and \eqref{190} are implied by
\cite[Theorem 5.2]{HLS}. The formulae \eqref{27} and \eqref{28}
can be derived as in the proof of the latter theorem;
cf.\ the computation of $\rho_{d-1,d}$ and
of $\rho_{0,d}$ in \cite{HLS}.

As in the computation of $\rho_{i,j}$ in \cite{HLS} (for $i=j=n-1$)
we obtain that
\begin{align}\label{29}
\rho(V_{g,n-1},V_{h,n-1})=A_0+A_1,
\end{align}
where
\begin{align*}
A_0:=\gamma^2&\iiiint e^{\gamma C_n(y-z)}\1\{y\in K^0_2,z\in K^0_1\}g(u)h(v)\\
&\times\Lambda_{n-1}(K_1;d(y,u))\,\Lambda_{n-1}(K_2;d(z,v))\,\Theta(dK_1)\,\BQ(dK_2)
\end{align*}
and
\begin{align*}
  A_1:=\gamma&\iiint e^{\gamma C_n(y-z)}g(u)h(v)\,
\Lambda_{n-1}(K;d(y,u))\,\Lambda_{n-1}(K;d(z,v))\,\BQ(dK).
\end{align*}
An easy calculation based on the covariance property of $\Lambda_{n-1}$ shows
that
\begin{align*}
A_0=\gamma^2 \int e^{\gamma C_n(x-z)}C_{n-1}(x-z;g)h(v)\,N_{n-1,n}(d(z,v,x)).
\end{align*}
As the number $A_1$ can be expressed directly as an integral
with respect to $N_{n-1,n-1}$, \eqref{31} follows.    \hfill \hbox{}\qed

\medskip

The following theorem establishes the existence of
asymptotic covariances for general geometric functionals.
Moreover, formula \eqref{sigmaphipsi} provides a tool
for their computation.

\begin{theorem}\label{thm:CovariancesGeneral}
Assume that \eqref{eqn:AssumptionSecondMoments} is satisfied and let 
$\psi$ and $\phi$ be geometric functionals. Then the limit
$$
\sigma(\psi,\phi):=\lim\limits_{r(W)\to\infty}\frac{\cov(\psi(Z\cap W),\phi(Z\cap W))}{V_n(W)}
$$
exists and is given by
\begin{align}\label{sigmaphipsi}
\sigma(\psi,\phi)=\varrho(\psi^*,\phi^*).
\end{align}
If \eqref{eqn:AssumptionThirdMoments} holds, there is a constant
$c_\Theta$, depending only on $\Theta$, such that, for
$W\in\cK^n$ with $r(W)\geq 1$,
\begin{equation}\label{eqn:ErrorCovarianceGeneral}
  \bigg| \frac{\cov(\psi(Z\cap W),\phi(Z\cap W))}{V_n(W)}-\sigma(\psi,\phi)\bigg|\leq
\frac{c_\Theta M(\psi) M(\phi)}{r(W)}.
\end{equation}
\end{theorem}

Theorem \ref{thm:CovariancesGeneral} is taken from \cite[Theorem
3.1]{HLS}. Its proof is involved and depends on the Fock space
representation~\cite{LaPe11} and several non-trivial integral-geometric inequalities for geometric
functionals. The inequality \eqref{eqn:ErrorCovarianceGeneral} allows
us to control the error if we approximate the exact covariance for a
given observation window by the asymptotic covariance. By evaluating
the left-hand side of \eqref{eqn:ErrorCovarianceGeneral} for the
volume one obtains a lower bound of order $1/r(W)$ (see
\cite[Proposition 3.8]{HLS}), which shows that the rate on the right-hand
side of \eqref{eqn:ErrorCovarianceGeneral} is optimal, in general.

Using Propositions \ref{prop1} and \ref{prop2}
in formula \eqref{sigmaphipsi}, we obtain the following result
for the asymptotic covariances
involving volume and surface content.

\begin{theorem}\label{thm:CovariancesVolumeSurfaceArea}
Assume that \eqref{eqn:AssumptionSecondMoments} holds and let $g,h:\S^{n-1}\to\R$ be measurable
and bounded. Then,
\begin{align*}
\sigma(V_n,V_n)&=(1-p)^2 \int (e^{\gamma C_n(x)} - 1) \, dx,\\
\sigma(V_{g,n-1},V_n) &= -(1-p)^2 \gamma \overpsi_{n-1}(g) \int (e^{\gamma C_n(x)} - 1) \, dx\\
& \qquad + (1-p)^2 \gamma \int g(u) e^{\gamma C_n(x-y)} \, N_{n-1,n}(d(x,u,y)).
\end{align*}
If, in addition, $\P(V_n(Z_0)>0)=1$, then
\begin{align*}
\sigma&(V_{g,n-1},V_{h,n-1}) = (1-p)^2\gamma^2 \overpsi_{n-1}(g)\overpsi_{n-1}(h)
\int (e^{\gamma C_n(x)}-1) \, dx\\
& + (1-p)^2 \gamma^2 \int e^{\gamma C_n(x-y)} h(u)C_{n-1}(x-y;g) \, N_{n-1,n}(d(y,u,x))\\
& -(1-p)^2 \gamma^2 \int e^{\gamma C_n(x-y)}
\big(g(u)\overpsi_{n-1}(h)+h(u)\overpsi_{n-1}(g)\big) \, N_{n-1,n}(d(y,u,x))\\
& + (1-p)^2 \gamma \int  e^{\gamma C_n(x-y)} g(u)h(v)\,N_{n-1,n-1}(d(x,u,y,v)).
\end{align*}
\end{theorem}

In the case $h=g\equiv 1$ the formula for 
$\sigma(V_{g,n-1},V_{h,n-1})$ simplifies to~\cite[Corollary 6.2]{HLS},
that is
\begin{align}\label{sigma11}\notag
\sigma (V_{n-1},V_{n-1}) =\; & (1-p)^2\gamma^2 v^2_{n-1}\int (e^{\gamma C_n(x)}-1) \, dx\\ \notag
& + (1-p)^2 \gamma^2 \int e^{\gamma C_n(x-y)}(C_{n-1}(x-y)-2v_1) \, N_{n-1,n}(d(y,u,x))\\
& + (1-p)^2 \gamma \int  e^{\gamma C_n(x-y)}\,N_{n-1,n-1}(d(x,u,y,v)),
\end{align}
where $C_{n-1}(x):=C_{n-1}(x;1)$ is defined by
Eq.~\eqref{eq:def_C_n_minus_1} with $g\equiv 1$.

In the planar case (treated in Subsection
\ref{subsec:CovariancesAnisotropic}) we will complement Theorem
\ref{thm:CovariancesVolumeSurfaceArea} with the asymptotic covariances
involving the Euler characteristic.  Integral representations of
asymptotic covariances of intrinsic volumes in general dimensions (with
respect to some special curvature based measures) can be found in
\cite[Sections 5 and 6]{HLS}.

Theorem \ref{thm:CovariancesGeneral} establishes the existence of an
asymptotic covariance matrix
$\Sigma=(\sigma(\psi_i,\psi_j))_{i,j=1,\hdots,m}$ for geometric
functionals $\psi_1,\hdots,\psi_m$. It is natural to ask whether this
matrix is positive definite. The next result (see~\cite[Theorem
4.1]{HLS}) gives sufficient, but presumably not necessary conditions for
positive definiteness.

\begin{theorem}\label{thm:PositiveDefiniteness}
  Let \eqref{eqn:AssumptionSecondMoments} be satisfied and assume that
$\P(V_n(Z_0)>0)>0$.
Let $\psi_0,\hdots,\psi_n$ be geometric functionals
  such that, for $i\in\{0,\hdots,n\}$, $\psi_i$ is homogeneous of
  degree~$i$ (that is, $\psi_i(\lambda K)=\lambda^i\psi_i(K)$ for 
  $\lambda >0$) and satisfies
$$
|\psi_i(K)|\geq \tilde{\beta}(\psi_i) r(K)^i, \quad K\in\cK^n,
$$
with a constant $\tilde{\beta}(\psi_i)$ only depending on
$\psi_i$. Then $\Sigma=(\sigma(\psi_i,\psi_j))_{0\leq i,j \leq n}$ is
positive definite.
\end{theorem}

Since the intrinsic volumes satisfy the assumptions of Theorem
\ref{thm:PositiveDefiniteness}, we obtain the following corollary.

\begin{corollary}
  Let \eqref{eqn:AssumptionSecondMoments} be satisfied and assume that
  the typical grain has nonempty interior with positive
  probability. Then the matrix
  $\Sigma=(\sigma(V_i,V_j))_{i,j=0,\hdots,n}$ is positive definite.
\end{corollary}

\subsection{Covariance formulae for planar Boolean
  models}\label{subsec:CovariancesAnisotropic}

In this section we consider the Boolean model in the planar
case $n=2$. For measurable and bounded $g:\S^1\to\R$
we consider the additive and measurable functional
\begin{align*}
V_{g,1}(K):=\Psi_1(K;g):=\int g(u)\,\Psi_1(K;du),\quad K\in\cR^n,
\end{align*}
see \eqref{Vgi}. We will compute the asymptotic covariances
between $V_0$  and the vector $(V_0,V_{g,1},V_2)$.

We define a function $\bar{h}:\S^1\to\R$ by
\begin{align}\label{2500}
\bar{h}(u):=\int h(K^*,u)\,\BQ(dK),\quad u\in\S^1,
\end{align}
where $K^*:=-K$ and $h(K^*,\cdot)$ is the support function of $K^*$. 
Indeed, if $K$ is a convex body containing the origin, then
the basic properties of $V_1$
together with the definition of the support function easily imply that
$0\le h(K^*,u)\le cV_1(K^*)=cV_1(K)$ for a constant $c>0$ that does only
depend on the dimension.  Therefore dominated convergence implies that
$\bar{h}$ is continuous and in particular bounded.  We also define
\begin{align}\label{2501}
v_{1,1}:=\overpsi_1(\bar{h})=\int \bar{h}(u)\, \overpsi_1(du)
=\iiint h(K^*,u)\,\overline{\Psi}_1(L;du)\,\BQ(dK)\,\BQ(dL).
\end{align}

\begin{theorem} Assume that \eqref{eqn:AssumptionSecondMoments}
and $\P(V_2(Z_0)>0)=1$ hold and let $g:\S^{1}\to\R$ be measurable
and bounded. Then
\begin{align}\label{45}\notag
\sigma(V_0,V_2)&=p(1-p)-(1-p)^2\gamma(1-\gamma v_{1,1})
\int \big(e^{\gamma C_2(x)}-1)\,dx\\
&\quad-2(1-p)^2\gamma^2\int \bar{h}(u)e^{\gamma
C_2(y-z)}\,N_{1,2}(d(y,u,z)),\allowdisplaybreaks\\
\label{46}\notag
\sigma(V_0,V_{g,1})&=(1-p)^2\gamma \overline{\Psi}_1(g)
+(1-p)^2\gamma^2\overline{\Psi}_1(g)(1-\gamma v_{1,1})\int\big(e^{\gamma C_2(x)}-1\big)\,dx\\
\notag
&\quad+(1-p)^2\int\big(\chi'(y-z)g(u)+2\gamma^3\overline{\Psi}_1(g)
e^{\gamma C_2(y-z)}\bar{h}(u)\big)\,N_{1,2}(d(z,u,y))\\
&\quad-2(1-p)^2\gamma^2\int e^{\gamma C_2(y-z)}
\bar{h}(u)g(v)\,N_{1,1}(d(y,u,z,v)),\allowdisplaybreaks\\
\label{47}\notag
\sigma(V_0,V_0)&=(1-2p)(1-p)\gamma +(1-p)(2p-3)v_{1,1}\gamma^2 \\\notag
&\quad+(1-p)^2\gamma^2(1-\gamma v_{1,1})^2
\int \big(e^{\gamma C_2(x)}-1\big)\,dx\\\notag
&\quad +(1-p)^2\int \bar{h}(u)\chi''(y-z)\,N_{1,2}(d(z,u,y))\\
&\quad +4(1-p)^2\gamma^3  \int e^{\gamma C_2(y-z)}\bar{h}(u)\bar{h}(v)\,N_{1,1}(d(y,u,z,v)),
\end{align}
where
\begin{align*}
\chi'(x)&:=e^{\gamma C_2(x)}
\big(\gamma^3(v_{1,1}-2  C_1(x;\bar{h}))-\gamma^2\big),\quad x\in\R^2,\\
  \chi''(x)&:=e^{\gamma C_2(x)}
\big(4\gamma^4(C_1(x;\bar{h})-v_{1,1})+4\gamma^3\big),\quad x\in\R^2.\\
\end{align*}
\end{theorem}
{\em Proof.} We wish to apply \eqref{sigmaphipsi}.
In view of Proposition \ref{prop1} we need to determine
$V^*_0$.
To do so we consider \eqref{11} for $j=0$ and $g\equiv 1$.
For the summation we distinguish
four cases. In the first two cases we have $m_i=0$ for exactly
one $i\in\{0,\hdots,k\}$ and either $m_0=0$ or $m_0=2$.
In the third and fourth case we have $m_i=m_r=1$ for exactly two
$i,r\in\{0,\hdots,k\}$ and either $m_0=0$ or $m_0=1$.
Accordingly we can write
\begin{align*}
\E V_0(Z\cap K_0)
=\sum^\infty_{k=1}\frac{(-1)^{k-1}\gamma^k}{k!}(a_{1,k}+a_{2,k}+a_{3,k}+a_{4,k}).
\end{align*}
The decomposability property \eqref{4264} (for $n=2$) 
and the symmetry of mixed functionals imply that
\begin{align*}
a_{1,k}=V_0(K_0)v^k_2,\quad a_{2,k}=V_2(K_0)kv^{k-1}_2.
\end{align*}
To treat $a_{3,k}$ and $a_{4,k}$ we use the decomposability property
\begin{align*}
  V^{(0)}_{1,1,2\dots,2}(K_0,\hdots,K_k)=V^{(0)}_{1,1}(K_0,K_1)V_2(K_2)\cdots V_2(K_k)
\end{align*}
and again the symmetry of mixed functionals (see \cite{Hug99}) to obtain that
\begin{align*}
a_{3,k}=\; &kv^{k-1}_2\int V^{(0)}_{1,1}(K_0,K)\,\BQ(dK),\\
a_{4,k}=\; &\frac{k(k-1)}{2}v^{k-2}_2V_2(K_0)\int V^{(0)}_{1,1}(K,L)\,\BQ^2(d(K,L)).
\end{align*}
It follows that
\begin{align}\label{13}\notag
  V^*_0(K_0)=\; &-(1-p)V_0(K_0)+(1-p)\gamma\int V^{(0)}_{1,1}(K_0,K)\,\BQ(dK)\\
&+(1-p)V_2(K_0)\bigg(\gamma-\frac{\gamma^2}{2}\int V^{(0)}_{1,1}(K,L)\,\BQ^2(d(K,L))\bigg)
\end{align}
or
\begin{align}\label{14}\notag
  V^*_0(K_0)=\; &-(1-p)V_0(K_0)+(1-p)\gamma\overline{V}_{1,1}(K_0)\\
&+(1-p)\Big(\gamma-\frac{\gamma^2}{2} w_{1,1}\Big)V_2(K_0),
\end{align}
where
\begin{align}
\overline{V}_{1,1}(K_0)&:=\int V^{(0)}_{1,1}(K_0,K)\,\BQ(dK),\\
w_{1,1}&:=\int V^{(0)}_{1,1}(K,L)\,\BQ^2(d(K,L)).
\end{align}

Using \eqref{sigmaphipsi} together with \eqref{14}
and Proposition \ref{prop1}, we obtain 
the following intermediate formulae
for the asymptotic covariances:
\begin{align}\label{15}\notag
\sigma(V_0,V_2)&=(1-p)^2\rho(V_0,V_2)-(1-p)^2\gamma\rho(\overline{V}_{1,1},V_2)\\
&\quad-(1-p)^2\Big(\gamma-\frac{\gamma^2}{2} w_{1,1}\Big)\rho(V_2,V_2),
\end{align}
\begin{align}\label{16}\notag
\sigma(V_0,V_{g,1})&=-(1-p)^2\gamma\overpsi_1(g)\rho(V_0,V_2)
+(1-p)^2\gamma^2\overpsi_1(g)\rho(\overline{V}_{1,1},V_2)\\ \notag
&\quad+(1-p)^2\overpsi_1(g)\Big(\gamma^2-\frac{\gamma^3}{2} w_{1,1}\Big)\rho(V_2,V_2)\\
\notag &\quad+(1-p)^2\rho(V_0,V_{g,1})
-(1-p)^2\gamma\rho(\overline{V}_{1,1},V_{g,1})\\
&\quad-(1-p)^2\Big(\gamma-\frac{\gamma^2}{2} w_{1,1}\Big)\rho(V_2,V_{g,1}),
\end{align}
and
\begin{align}\label{17}
\notag
\sigma(V_0,V_0)
&=(1-p)^2\rho(V_0,V_0)-2(1-p)^2\gamma\rho(V_0,\overline{V}_{1,1})\\\notag
&\quad+(1-p)^2\gamma^2\rho(\overline{V}_{1,1},\overline{V}_{1,1})
-2(1-p)^2\Big(\gamma-\frac{\gamma^2}{2} w_{1,1}\Big)\rho(V_0,V_2)\\\notag
&\quad+2(1-p)^2\Big(\gamma^2-\frac{\gamma^3}{2} w_{1,1}\Big)\rho(\overline{V}_{1,1},V_2)\\
&\quad 
+(1-p)^2\Big(\gamma-\frac{\gamma^2}{2} w_{1,1}\Big)^2\rho(V_2,V_2).
\end{align}

At this stage we can use the formula
\begin{align*}
  V^{(0)}_{1,1}(K,L)=2\int h(L^*,u)\,\Psi_1(K;du),\quad K,L\in\cK,
\end{align*}
(use (6.25) and (14.21) in \cite{SchneiderWeil} as well as
$S_1=2\Psi_1$) implying that
\begin{align}\label{211}
\overline{V}_{1,1}(K)&=2\int \bar{h}(u)\,\Psi_1(K;du)=2 V_{\bar{h},1}(K),\\
\label{212}
w_{1,1}&= \int\overline{V}_{1,1}(K)\,\BQ(dK)=2\overpsi_1(\bar{h})=2v_{1,1}.
\end{align}

Theorem 5.2 in \cite{HLS} implies that
\begin{align}\label{191}
\rho(V_0,V_0)=e^{\gamma v_2}\Big(\gamma+\frac{\gamma^2 v'_{1,1}}{2}\Big),
\end{align}
where
\begin{align}\label{192}
v'_{1,1}:=\int \Phi_0(K_1\cap (K_2+x);\partial{K}_1\cap(\partial{K}_2+x))\,dx\,\BQ^2(d(K_1,K_2)).
\end{align}
It follows from \cite[Theorem 6.4.1]{SchneiderWeil} (together with the decomposability property and the fact
that the boundary of a convex body has vanishing volume) that
$$
\int \Phi_0(K_1\cap (K_2+x);\partial{K}_1\cap(\partial{K}_2+x))\,dx
=\Phi^{(0)}_{1,1}(K_1,K_2;\partial K_1\times\partial K_2),
$$
where $\Phi^{(0)}_{1,1}(K_1,K_2;\cdot)$ is a mixed functional. 
Since $\Phi^{(0)}_{1,1}(K_1,K_2;\cdot)$ is concentrated on $\partial K_1\times\partial K_2$ by 
\cite[Theorem 6.4.1 (b)]{SchneiderWeil}, we have
\begin{align*}
\Phi^{(0)}_{1,1}(K_1,K_2;\partial K_1\times\partial K_2)=\Phi^{(0)}_{1,1}(K_1,K_2;\R^2\times\R^2)
=V^{(0)}_{1,1}(K_1,K_2).
\end{align*}
Therefore $v'_{1,1}=w_{1,1}=2v_{1,1}$ and
\begin{align}\label{1911}
  \rho(V_0,V_0)=(1-p)^{-1}(\gamma+\gamma^2 v_{1,1}).
\end{align}

Now we can insert \eqref{211} and \eqref{212}
as well as \eqref{1911} and the formulae of Proposition \ref{prop2} into
\eqref{15}--\eqref{17} to obtain the assertions. From \eqref{15} we get
\begin{align*} 
\sigma(V_0,V_2)&=(1-p)^2\rho(V_0,V_2)-(1-p)^2\gamma 2\rho(V_{\bar h,1},V_2)\\
&\qquad -(1-p)^2(\gamma-\gamma^2 v_{1,1})\rho(V_2,V_2)
\end{align*}
so that \eqref{45} follows from \eqref{190}, \eqref{27} and \eqref{VnVn}. 

Next we deduce from \eqref{16} that
\begin{align*}
\sigma(V_0,V_{g,1})&=-p(1-p)\gamma\overpsi_1(g)\\
&\quad+2(1-p)^2\gamma^3\overpsi_1(g)\int \bar{h}(u)e^{\gamma C_2(y-z)}\,N_{1,2}(d(y,u,z))\\
&\quad+(1-p)^2\overpsi_1(g)(\gamma^2-\gamma^3 v_{1,1})
\int \big(e^{\gamma C_2(x)}-1)\,dx\\
&\quad+(1-p)^2\gamma e^{\gamma v_2}\overpsi_1(g)\\
&\quad-2(1-p)^2\gamma^3 \int e^{\gamma C_2(y-z)}C_1(y-z;\bar{h})g(v)\,
N_{1,2}(d(z,v,y))\\
&\quad-2(1-p)^2\gamma^2 \int e^{\gamma C_2(y-z)}\bar{h}(u)g(v)\,N_{1,1}(d(y,u,z,v))\\
&\quad-(1-p)^2\Big(\gamma^2-\gamma^3 v_{1,1}\Big)\int g(u)e^{\gamma C_2(y-z)}\,N_{1,2}(d(y,u,z)).
\end{align*}
Equation \eqref{46} follows upon some simplification and rearrangement.

From  \eqref{17} we obtain that
\begin{align*}
\sigma(V_0,V_0)&=(1-p)\big(\gamma+\gamma^2 v_{1,1}\big)-4(1-p)\gamma^2v_{1,1}\\
 &\quad+4(1-p)^2\gamma^4 \int e^{\gamma C_2(y-z)}C_1(y-z;\bar{h})\bar{h}(v)\,
N_{1,2}(d(z,v,y))\\
&\quad+4(1-p)^2\gamma^3 \int e^{\gamma C_2(y-z)}\bar{h}(u)\bar{h}(v)N_{1,1}(d(y,u,z,v))\\
&\quad-2(1-p)^2\big(\gamma-\gamma^2 v_{1,1}\big)((1-p)^{-1}-1)\\
&\quad+4(1-p)^2(\gamma^3-\gamma^4 v_{1,1})
\int \bar{h}(u)e^{\gamma C_2(y-z)}\,N_{1,2}(d(y,u,z))\\
&\quad+(1-p)^2\big(\gamma-\gamma^2 v_{1,1}\big)^2\int \big(e^{\gamma C_2(x)}-1)\,dx.
\end{align*}
Equation \eqref{47} now follows from an easy calculation.    \hfill \hbox{}\qed

\medskip

In the isotropic case, $\bar{h}:=\bar{h}(u)$ does not depend on
$u\in\S^1$. By \cite[(14.21)]{SchneiderWeil}, we have for $L\in\cK^n$ 
\begin{align*}
V_1(L)=\int h(L,u)\,\Psi_1(B^2;du)=\frac{1}{2}\int h(L,u)\,\cH^1(du),
\end{align*}
so that
$$
v_1=\frac{1}{2}\iint h(L,u)\,\cH^1(du)\,\BQ(dL)=\pi\bar{h}.
$$
Further
$$
v_{1,1}=\iint \bar{h}(u)\,\,\Psi_1(K;du)\,\BQ(dK)=\bar{h}v_1.
$$
Hence
\begin{align}\label{33}
\bar{h}=\frac{v_1}{\pi},\quad v_{1,1}=\frac{v^2_1}{\pi}.
\end{align}
Inserting \eqref{33} into \eqref{45}--\eqref{47} yields Corollary 6.3 in \cite{HLS}.

\subsection{The Boolean model on the torus}
\label{subsec:torus}

\begin{figure}[b]
 \centering
 \includegraphics[width=0.5\textwidth]{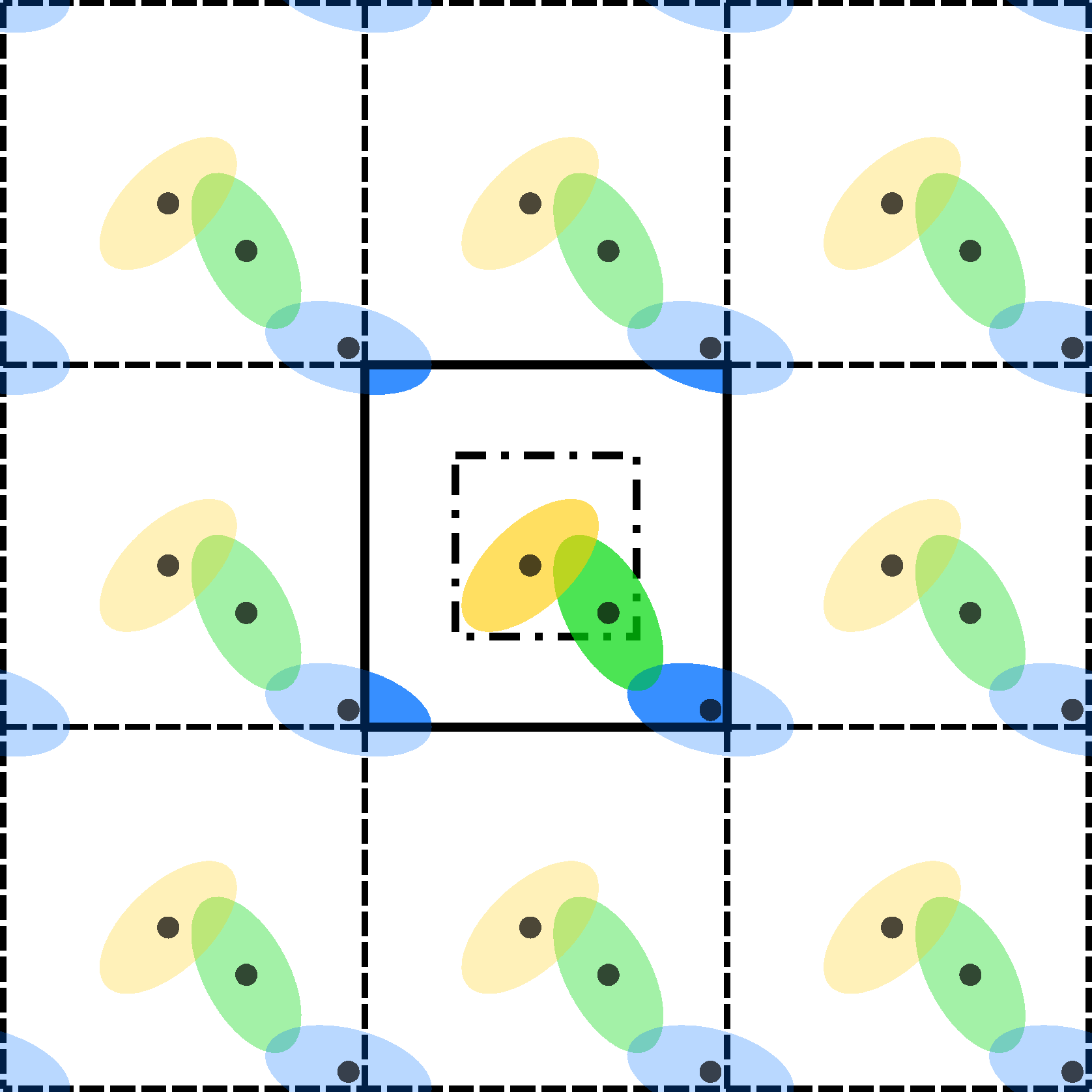}
 \caption{The Boolean model with periodic boundary conditions: we
 consider grains with centers in $[-1/2,1/2]^n$ (square with solid line) and all
 their translations by $\mathbb{Z}^n$ valued vectors (in squares with
 dashed lines). The Boolean model with periodic boundary conditions is
 obtained by taking the union of all the grains and restricting to the
 square with the solid line.}
 \label{fig:pbc}
\end{figure}

One obtains the $n$-dimensional (unit) torus $T^n$ by identifying
opposite sides of the boundary of $[-1/2,1/2]^n$. As in $\R^n$ one can
consider a translation invariant Poisson process of grains
on the torus (with intensity measure $\Lambda$, grain distribution
$\mathbf{Q}$ and intensity $\gamma$) and consider the resulting Boolean
  model $Z_{T^n}$.
The Boolean model on the torus $T^n$ can be constructed in the following
way from a random closed set in $\R^n$ (see Fig.~\ref{fig:pbc}). We
start with a homogeneous Poisson process in $[-1/2,1/2]^n$ and put
around each point an independent copy of the typical grain. For each
grain, we also place all translates by vectors $v\in\mathbb{Z}^n$
and take the union of all resulting grains. Finally, we restrict this random
closed set to $[-1/2,1/2]^n$ and identify opposite boundaries. This
setting is also denoted as periodic boundary conditions.

For a geometric functional $\psi:\mathcal{R}^n\to\R$ we can define
$\psi(Z_{T^n})$ in the following way. For a set $K\subset T^n$ whose
embedding $K_{\R^n}=\{x\in\R^n: x\in K\}$ into $\R^n$ is a convex body
and is contained in $(-1/2,1/2)^n$ we put $\psi(K)=\psi(K_{\R^n})$. By
further requiring that $\psi$ is translation-invariant and additive on
$T^n$, this gives us $\psi(Z_{T^n})$.

By computing the Fock space representation of $\psi(Z_{T^n})$ and
$\phi(Z_{T^n})$ for geometric functionals $\psi,\phi:\mathcal{R}^n\to\R$ as
in \cite[Section 3]{HLS} for a Boolean model in $\R^n$, one obtains that
\begin{align*}
& \cov(\psi(Z_{T^n}),\phi(Z_{T^n}))\\
& =\sum_{n=1}^\infty \frac{\gamma}{n!} \iint  (\E \psi(Z_{T^n}\cap K_1\cap \hdots\cap K_n)-\psi(K_1\cap \hdots\cap K_n))\\
 & \qquad \times (\E \phi(Z_{T^n}\cap K_1\cap \hdots\cap K_n)-\phi(K_1\cap \hdots\cap K_n)) \, \Lambda^{n-1}(d(K_2,\hdots,K_n)) \, \BQ(dK_1).
\end{align*}
Now let us assume that the grain distribution $\BQ$ is such that the
typical grain $Z_0$ is almost surely contained in $[-1/4,1/4]^n$, which
is depicted by the dot-dash line in Fig.~\ref{fig:pbc}. In this case the
intersection of two grains is always convex and the intersections on the
right-hand side of the covariance formula are the same as for a Boolean
model in $\R^n$ with grain distribution $\BQ$ and intensity $\gamma$.
Thus, it follows from the above definition of $\psi$ and $\phi$ of a
subset of the torus whose embedding into $\R^n$ is a convex body and a
subset of $(-1/2,1/2)^n$ and the additivity that
$$
\cov(\psi(Z_{T^n}),\phi(Z_{T^n}))=\rho(\psi^*,\phi^*).
$$
In other words, if the typical grain is sufficiently bounded, the exact
covariances for the Boolean model on the torus coincide with the
asymptotic covariances for the corresponding Boolean model in $\R^n$.
This provides a way to compute estimates for the asymptotic covariances
via simulations on the torus.

\section{Central limit theorems}\label{sec:CLTs}

In this section we consider the asymptotic behaviour of the
distributions of geometric functionals or of vectors of geometric
functionals for growing observation window.
Recall that a sequence of $m$-dimensional random vectors
$(Y_i)_{i\in\N}$ converges in distribution to an $m$-dimensional
random vector $Y$ if
$$
\lim_{i\to\infty} \P(Y_i\leq x) = \P(Y\leq x)
$$
for all $x\in\R^m$ such that $y\mapsto \P(Y\leq y)$ is continuous at $x$. 
(Here the relation $\le$ is to be understood componentwise).
In this case we write $Y_i\overset{d}{\longrightarrow} Y$ (as $i\to\infty$).
We are not only interested in the convergence in distribution 
but also in error bounds.
In order to measure the distance between
 the distributions of two $m$-dimensional random vectors $Z_1,Z_2$,
we use the $d_3$-metric which
is given by
$$
d_3(Z_1,Z_2)=\sup_{h\in \mathcal{H}_m} |\E h(Z_1)-\E h(Z_2)|,
$$
where $\mathcal{H}_m$ is the set of all $C^3$-functions $h:\R^m\to\R$
such that the absolute values of the second and the third partial
derivatives are bounded by one. For two random variables $Z_1,Z_2$ we
consider the Wasserstein distance
$$
d_W(Z_1,Z_2) = \sup_{h\in\operatorname{Lip}(1)} |\E h(Z_1) - \E h(Z_2)|,
$$
where $\operatorname{Lip}(1)$ is the set of all functions $h:\R\to\R$
whose Lipschitz constant is at most one. Note that convergence in the
$d_3$-distance or in the Wasserstein distance implies convergence in
distribution.

For the quantitative bounds we assume that there is a constant
$\varepsilon>0$ such that
\begin{equation}\label{eqn:AssumptionMoments3PlusEpsilon}
\E V_i(Z_0)^{3+\varepsilon}<\infty, \quad i=0,\hdots,n.
\end{equation}

We begin with a multivariate central limit theorem for a vector of geometric functionals.

\begin{theorem}\label{thm:MultivariateCase}
  Assume that \eqref{eqn:AssumptionSecondMoments} is satisfied, let
  $\Psi:=(\psi_1,\hdots,\psi_m)$ for geometric functionals
  $\psi_1,\hdots,\psi_m$, and let $N_\Sigma$ be an $m$-dimensional
  centred Gaussian random vector with covariance matrix
  $\Sigma=(\sigma(\psi_i,\psi_j))_{i,j=1,\hdots,m}$. Then
$$
\frac{1}{\sqrt{V_n(W)}} \big(\Psi(Z\cap W)- \E\Psi(Z\cap W) \big) 
\overset{d}{\longrightarrow}
N_\Sigma \quad \text{ as } \quad r(W)\to\infty.
$$
If \eqref{eqn:AssumptionMoments3PlusEpsilon} holds, there is a
constant $C_{\psi_1,\hdots,\psi_m}$ depending on
$\psi_1,\hdots,\psi_m$, $\Theta$ and $\varepsilon$ such that
$$
d_3\bigg(\frac{1}{\sqrt{V_n(W)}} \big(\Psi(Z\cap W)- \E\Psi(Z\cap W)
\big) ,N_\Sigma\bigg) \leq
\frac{C_{\psi_1,\hdots\psi_m}}{r(W)^{\min\{\varepsilon n/2,1\}}}
$$
for $W\in\cK^n$ with $r(W)\geq 1$.
\end{theorem}

This result was proved in \cite[Theorem 9.1]{HLS}
by using the Stein-Malliavin method and a truncation argument.

As tensors can be interpreted as vectors, we can define
convergence of tensor valued random elements and their $d_3$-distance
via convergence and $d_3$-distance for random vectors.  Since the
components of $\Phi_m^{0,s}$ are geometric functionals, Theorem
\ref{thm:MultivariateCase} can be applied to the translation invariant
Minkowski tensors.

\begin{corollary}\label{cor:CLTMinkowskiTensors}
  Assume that \eqref{eqn:AssumptionSecondMoments} holds, let $s\in\N$
  and $m\in\{0,\hdots,n-1\}$ and let $N^{0,s}_m$ be a random element
  in $\mathbb{T}^s$ such that each component is a centred Gaussian
  random variable and
$$
\cov((N)_{i_1,\hdots,i_s},(N)_{j_1,\hdots,j_s})=
\sigma((\Phi_m^{0,s})_{i_1,\hdots,i_s},(\Phi_m^{0,s})_{j_1,\hdots,j_s})
$$
for $i_1,\hdots,i_s,j_s,\hdots,j_s\in\{1,\hdots,n\}$. Then
$$
\frac{1}{\sqrt{V_n(W)}} (\Phi_m^{0,s}(Z\cap W)-\E\Phi_m^{0,s}(Z\cap W) ) 
\overset{d}{\longrightarrow} N^{0,s}_m \quad \text{ as } \quad r(W)\to\infty.
$$
If \eqref{eqn:AssumptionMoments3PlusEpsilon} holds, there is a
constant $C_{s,m}$ depending on $s$, $m$, $\Theta$ and $\varepsilon$
such that
$$
d_3\bigg(\frac{1}{\sqrt{V_n(W)}} (\Phi_m^{0,s}(Z\cap
W)-\E\Phi_m^{0,s}(Z\cap W) ) ,N_m^{0,s}\bigg) \leq
\frac{C_{s,m}}{r(W)^{\min\{\varepsilon n/2,1\}}}
$$
for $W\in\cK^n$ with $r(W)\geq 1$.
\end{corollary}

In the multivariate case we assume translation invariance of the
geometric functionals in order to ensure the existence of an
asymptotic covariance matrix. In the univariate case this is not
required since one can standardize be dividing by the standard
deviation. For this reason, we can drop the assumption of translation
invariance in the following univariate central limit theorem, which
is taken from \cite[Theorem 9.3]{HLS}.

\begin{theorem}\label{thm:UnivariateCase}
  Let \eqref{eqn:AssumptionSecondMoments} be satisfied, let
  $\psi:\mathcal{R}^n \to \R$ be measurable, additive and locally
  bounded, assume that there are constants $r_0 \geq 1$ and $\sigma_0>0$ such that
$$
\frac{\var \psi(Z\cap W)}{V_n(W)} \geq \sigma_0
$$
for $W\in\cK^n$ with $r(W)\geq r_0$ and let $N$ be a standard Gaussian random variable. Then
$$
\frac{\psi(Z\cap W)-\E \psi(Z\cap W)}{\sqrt{\var \psi(Z\cap W)}} 
\overset{d}{\longrightarrow} N \quad \text{ as } \quad r(W)\to\infty.
$$
If, additionally, \eqref{eqn:AssumptionMoments3PlusEpsilon} is
satisfied, there is a constant $c_\psi$ depending on $\psi$, $\Theta$,
$r_0$, $\sigma_0$, and $\varepsilon$ such that
$$
d_W\bigg(\frac{\psi(Z\cap W)-\E \psi(Z\cap W)}{\sqrt{\var \psi(Z\cap
    W)}},N\bigg)\leq \frac{c_\psi}{V_n(W)^{\min\{\varepsilon/2,1/2\}}}
$$
for $W\in\cK^n$ with $r(W)\geq r_0$.
\end{theorem}

The results presented in this section generalize previous findings in
\cite{BaryshnikovYukich,Baddeley1980,Heinrich2005,
HeinrichMolchanov1999, Mase1982, Molchanov1995, Penrose2007}, which only
deal with volume, surface area or closely related functionals.

\section{Boolean model of aligned
  rectangles}\label{sec:rectangles}

In this section we assume that $n=2$ and that the typical grain $Z_0$
is a deterministic rectangle of the form
$$
K:=\Big[-\frac{a}{2}e_1,\frac{a}{2}e_1\Big]+\Big[-\frac{b}{2}e_2,\frac{b}{2}e_2\Big]
= \Big[-\frac{a}{2},\frac{a}{2}\Big]\times\Big[-\frac{b}{2},\frac{b}{2}\Big]
$$
for some fixed $a,b>0$, where $e_1:=(1,0)$ and $e_2:=(0,1)$.
Then $v_2=ab$ and $v_1=a+b$.

\begin{figure}[t]
  \centering
  \subfigure[][]{\includegraphics[width=0.482\textwidth]{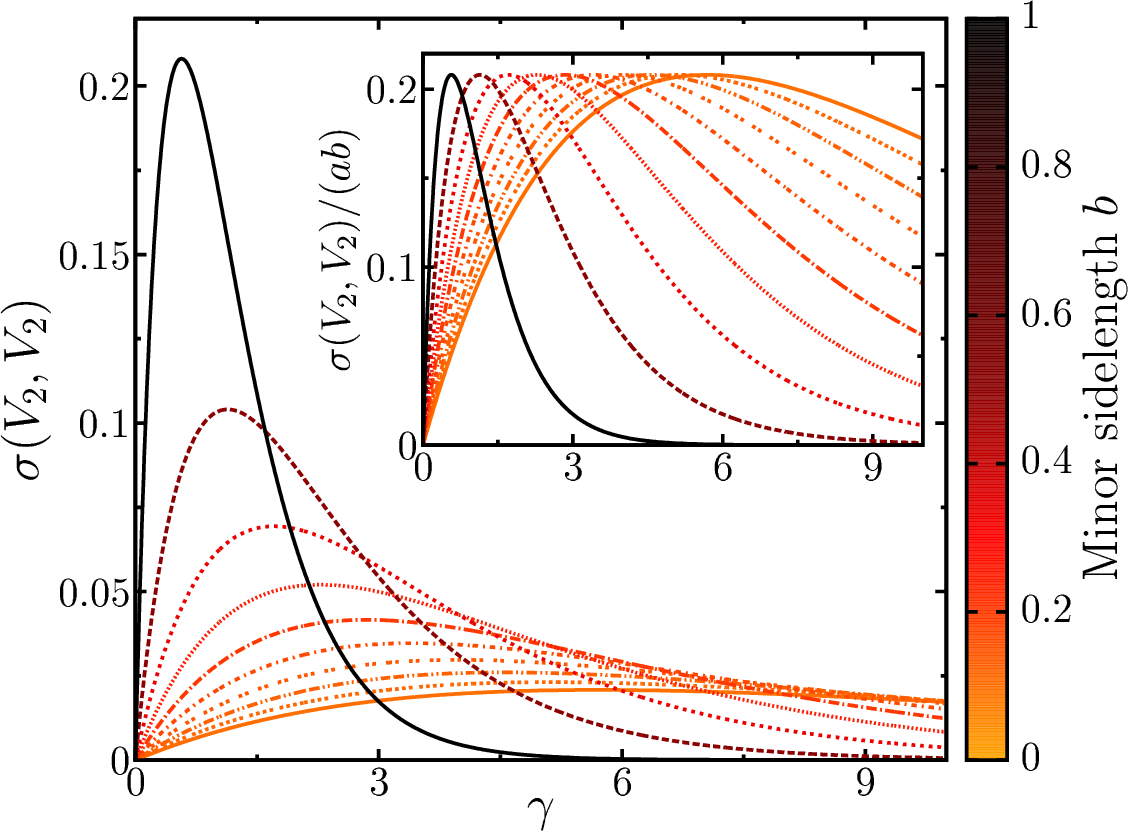}
  \label{var2colorall} }%
  \hfill%
  \subfigure[][]{\includegraphics[width=0.482\textwidth]{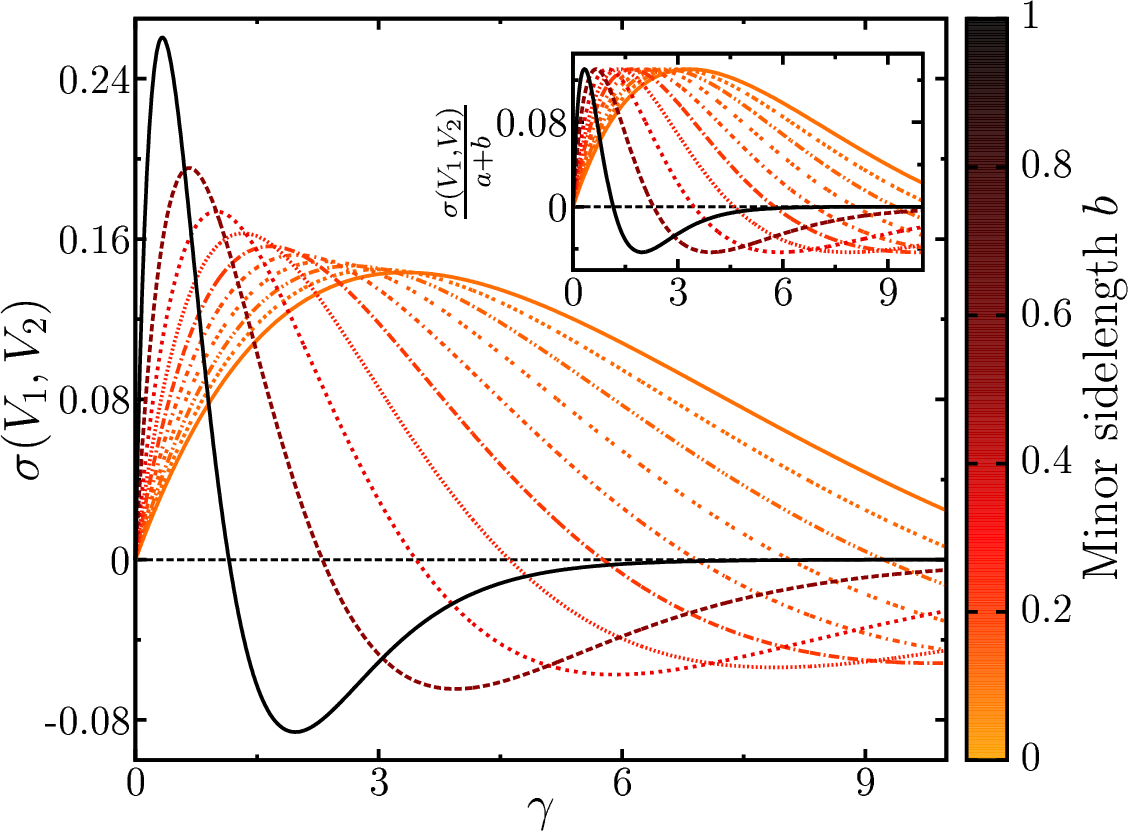}
  \label{cov12colorall} }\\
  \subfigure[][]{\includegraphics[width=0.482\textwidth]{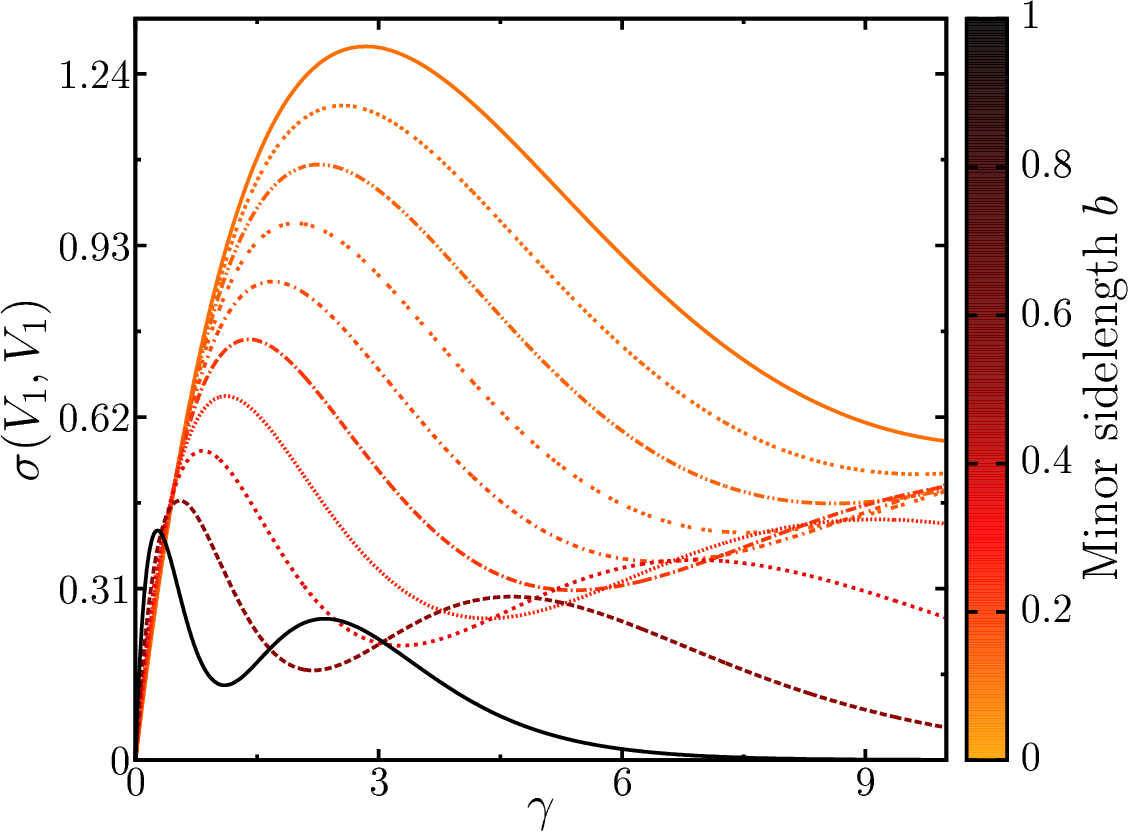}
  \label{var1colorall} }%
  \hfill%
  \subfigure[][]{\includegraphics[width=0.482\textwidth]{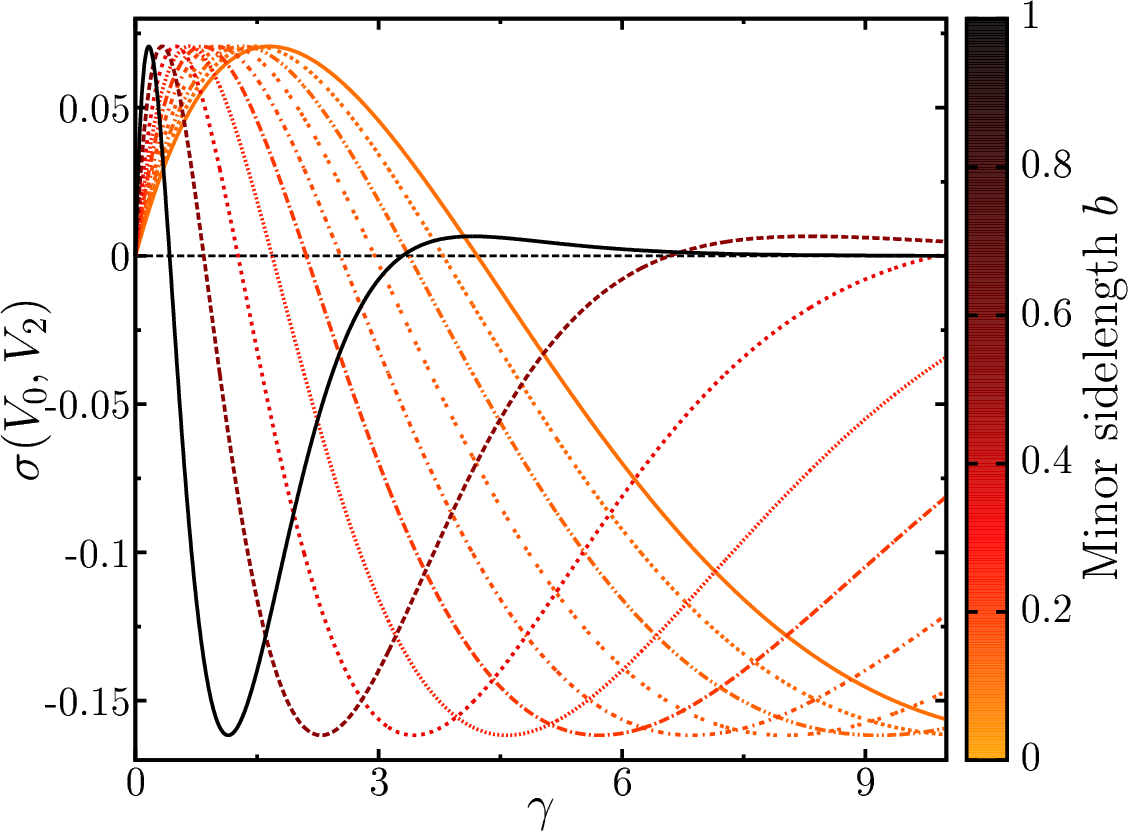}
  \label{cov02colorall} }\\
  \subfigure[][]{\includegraphics[width=0.482\textwidth]{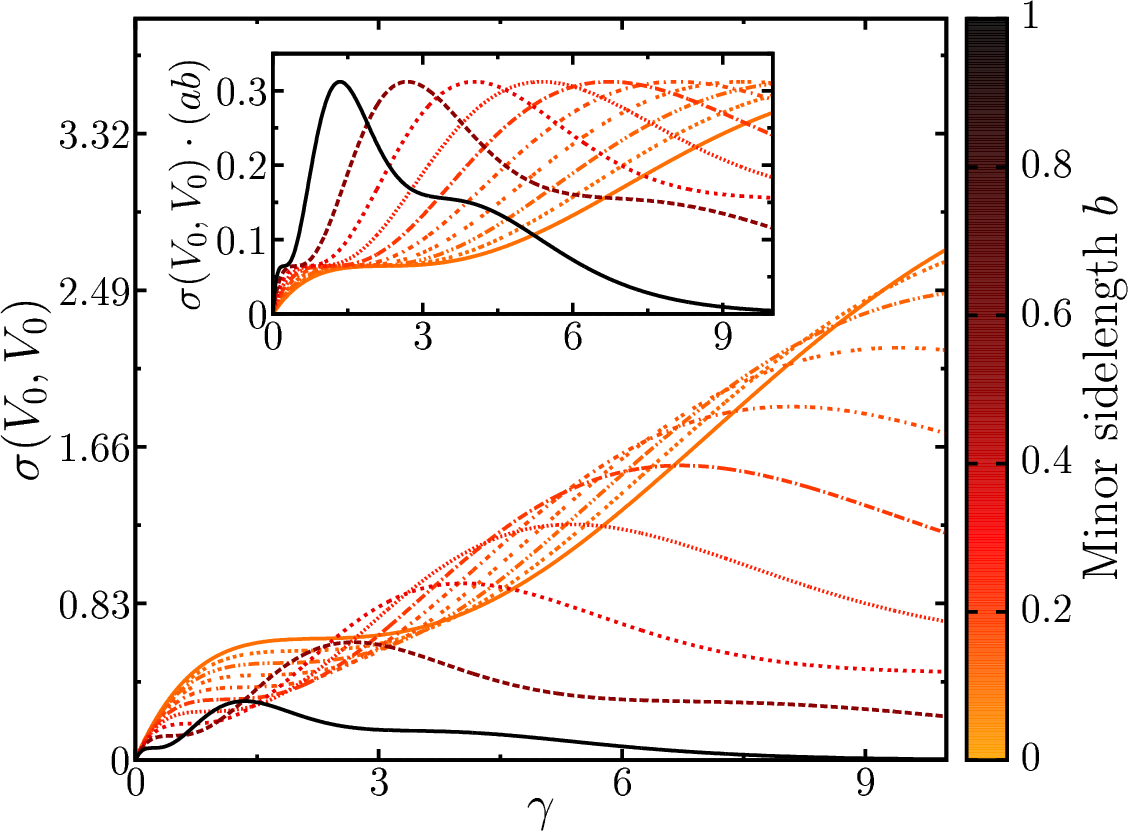}
  \label{var0colorall} }%
  \hfill%
  \subfigure[][]{\includegraphics[width=0.482\textwidth]{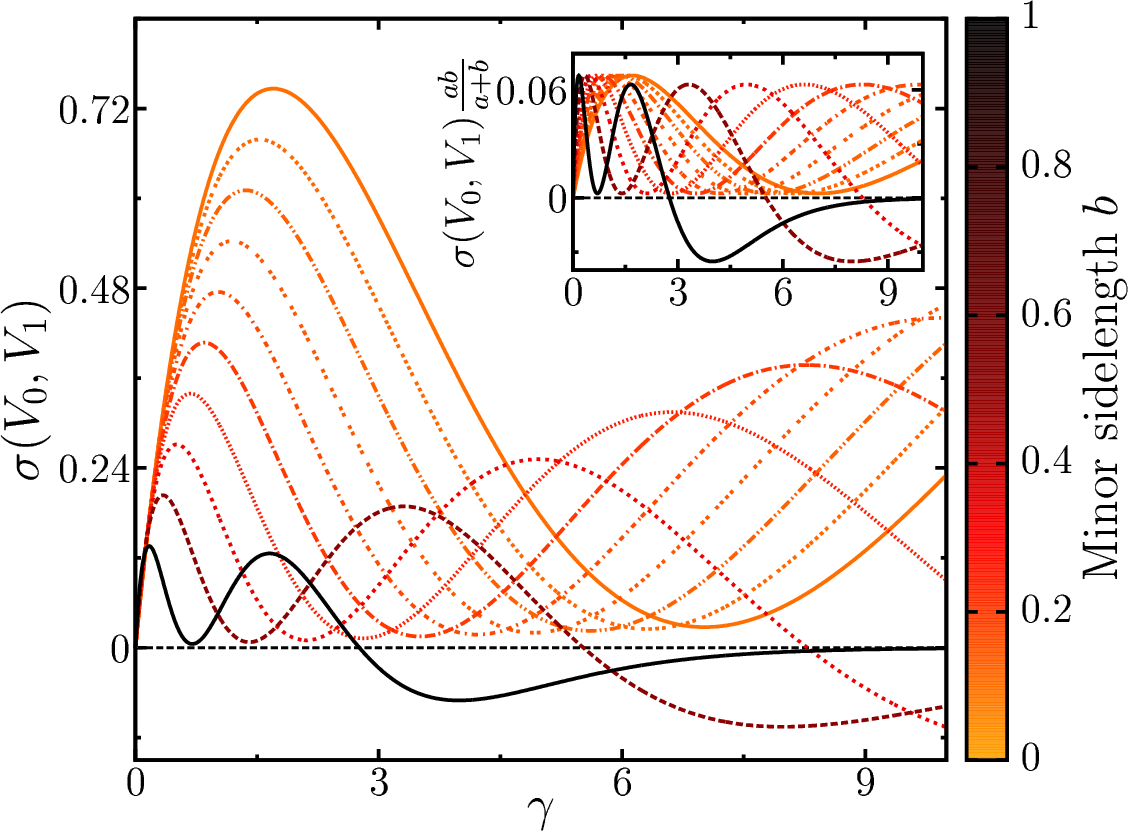}
  \label{cov01colorall} }
  \caption{Asymptotic covariances $\sigma(V_i,V_j)$ as a function of the
    intensity $\gamma$ for Boolean models of aligned rectangles with
    varying aspect ratio $b/a$; we choose $a=1$, hence $b\in(0,1]$; see Eqs.~\protect\eqref{rect2},
    \protect\eqref{rect5}, \protect\eqref{4522}, \protect\eqref{4527}, \protect\eqref{sigma01}, and
    \protect\eqref{rect6789}. The insets show covariances that are rescaled by
    suitable functions of the side lengths $a$ and $b$ of a single rectangle so that they only
    depend on $\gamma v_2$ but not on the aspect ratio, which also holds
    for $\sigma(V_0,V_2)$ in (d) without rescaling.}
   \label{fig_cov_numerical_integration}
\end{figure}

For any $x=(x_1,x_2)\in\R^2$ we have
$$
C_2(x)=V_2(K\cap (K+x))=\1\{|x_1|\le a,|x_2|\le b\}(a-|x_1|)(b-|x_2|).
$$
A change of variables and a symmetry argument imply that
\begin{align}\label{rect1}
  \int \big(e^{\gamma C_2(x)}-1\big)\,dx=4v_2H(\gamma v_2),
\end{align}
where the function $H:[0,\infty)\to [0,\infty)$ is defined by
\begin{align}\label{I11}
  H(r):=\int^1_0\int^1_0 \big(e^{rst}-1\big)\,ds\,dt
=\sum^\infty_{k=1}\frac{r^k}{k!(k+1)^2},\quad r\ge 0.
\end{align}
Hence we obtain from Theorem \ref{thm:CovariancesVolumeSurfaceArea}
that
\begin{empheq}[box=\fbox]{align}
  \label{rect2}
  \sigma(V_2,V_2)=4(1-p)^2v_2H(\gamma v_2),
\end{empheq}
where we recall that $p=1-e^{-\gamma v_2}$. The variance is visualized in
Fig.~\ref{var2colorall}.

At this stage it is convenient to complement the definition \eqref{I11}
with the following easy to check formulae:
\begin{align}\label{I12}
\int^1_0\int^1_0 e^{r st}s\,ds\,dt
&=\frac{1}{r^2}e^{r}-\frac{1}{r^2}-\frac{1}{r},\allowdisplaybreaks\\
\label{I13}
\int^1_0\int^1_0 e^{r
st}s^2\,ds\,dt&=\frac{1}{r^2}e^r-\frac{1}{r^3}e^r+\frac{1}{r^3}-\frac{1}{2r},\allowdisplaybreaks\\
\label{I14}
\int^1_0\int^1_0 e^{r st}st\,ds\,dt&=\frac{1}{r^2}e^r-\frac{1}{r^2}-\frac{1}{r}H(r)-\frac{1}{r}.
\end{align}
A consequence is
\begin{align} \label{I15}
\int^1_0\int^1_0 e^{r st}(st+s^2)\,ds\,dt&=\frac{2}{r^2}e^r-\frac{1}{r^3}e^r+\frac{1}{r^3}
-\frac{1}{r^2}-\frac{3}{2r}-\frac{1}{r}H(r).
\end{align}

Now we use Theorem \ref{thm:CovariancesVolumeSurfaceArea}
(for $n=2$ and $g\equiv 1$) to compute $\sigma(V_1,V_2)$.
For any measurable and even functions $f:\R^2\to[0,\infty)$
and $\tilde{f}:\S^1\to[0,\infty)$ we have \begin{align*} \int
    f(y-z)\tilde{f}(u) \, N_{1,2}(d(y,u,z)) =a^+_1+a^-_1+a^+_2+a^-_2,
  \end{align*} where \begin{align*} a^\pm_i:=\frac{1}{2}\iint \1\{y\in
    A^\pm_i,z\in K\}f(y-z)\tilde{f}(e_i)\,\cH^1(dy)\,dz,\quad
    i\in\{1,2\} \end{align*} and $A^\pm_1:=\{(x_1,x_2)\in K:x_1=\pm
  a/2\}$, $A^\pm_2:=\{(x_1,x_2)\in K:x_2=\pm b/2\}$.  By Fubini and a
  change of variables \begin{align}\label{rect9}
    a^\pm_1=\frac{\tilde{f}(e_1)}{2}\iint \1\{y\in A^\pm_1,y+z\in
    K\}f(z)\,\cH^1(dy)\,dz.  \end{align} For any $z=(z_1,z_2)\in K$ with
  $-a\le z_1\le 0$ and $z_2\ge 0$ we have \begin{align*} \int \1\{y\in
    A_1,y+z\in K\}\,\cH^1(dy)=\cH^1([-b/2,b/2-z_2])=b-z_2=b-|z_2|.
  \end{align*} For $-a\le z_1\le 0$ and $z_2\le 0$ this integral takes
  the same value.  Since the set of all $z$ with $z_1\notin[-a,0]$ or
  $|z_2|>b$ does not contribute to $a^+$ while the set of all $z$ with
  $z_1\notin[0,a]$ or $|z_2|>b$ does not contribute to $a^-$ it follows
  that \begin{align}\label{rect46}
    a^+_1+a^-_1=\; &\frac{\tilde{f}(e_1)}{2}\int \1\{|z_1|\le a,|z_2|\le
    b\}f(z_1,z_2)(b-|z_2|)\,d(z_1,z_2).  \end{align} Using a similar
  result for $b^+_1+b^-_1$ gives \begin{align}\label{rect52}\notag \int
    f(y-z)\tilde{f}(u)& \, N_{1,2}(d(y,u,z))\\ \notag
    =\; &\frac{\tilde{f}(e_1)}{2}\int \1\{|z_1|\le a,|z_2|\le
    b\}f(z_1,z_2)(b-|z_2|)\,d(z_1,z_2)\\ &+\frac{\tilde{f}(e_2)}{2}\int
    \1\{|z_1|\le a,|z_2|\le b\}f(z_1,z_2)(a-|z_1|)\,d(z_1,z_2).
  \end{align}
Inserting here $f(z):=e^{\gamma C_2(z)}$ and using a
change of variables gives
\begin{align*}
  \int &e^{\gamma C_2(y-z)}\tilde{f}(u) \, N_{1,2}(d(y,u,z))\\
    &= 2ab^2\tilde{f}(e_1)\int^1_0\int^1_0 e^{\gamma ab
    y_1y_2}y_2\,dy_1\,dy_2 +2a^2b\tilde{f}(e_2)\int^1_0\int^1_0
    e^{\gamma ab y_1y_2}y_1\,dy_1\,dy_2.
\end{align*}
From \eqref{I12} we obtain that
\begin{align}\label{rec652}\notag
  \int &e^{\gamma C_2(y-z)}\tilde{f}(u) \, N_{1,2}(d(y,u,z))\\
    &=\tilde{f}(e_1)\bigg(\frac{2}{\gamma^2a}e^{\gamma
    ab}-\frac{2}{\gamma^2a}-\frac{2b}{\gamma}\bigg)
    +\tilde{f}(e_2)\bigg(\frac{2}{\gamma^2b}e^{\gamma
    ab}-\frac{2}{\gamma^2b}-\frac{2a}{\gamma}\bigg).
\end{align}
In the case $\tilde{f}(e_1)=\tilde{f}(e_2)=1$ this yields
  \begin{align}\label{rect71} \int e^{\gamma C_2(x-y)} \,
    N_{1,2}(d(x,u,y)) =2v_1\bigg(\frac{1}{\gamma^2v_2}e^{\gamma
    v_2}-\frac{1}{\gamma^2v_2}-\frac{1}{\gamma}\bigg).  \end{align}
  Inserting this result together with \eqref{rect1} into the formula of
  Theorem \ref{thm:CovariancesVolumeSurfaceArea} yields
\begin{empheq}[box=\fbox]{align}
  \label{rect5} 
  \sigma(V_1,V_2)= 2(1-p)^2 v_1\left[\frac{1}{\gamma v_2}(e^{\gamma
    v_2}-1)-1-2 \gamma v_2H(\gamma
  v_2)\right],
\end{empheq}
which is visualized in Fig.~\ref{cov12colorall}.

Next we use \eqref{45} to compute $\sigma(V_0,V_2)$, starting with the
observation $$ h(K,-e_1)=h(K,e_1)=\frac{a}{2},\qquad
h(K,-e_2)=h(K,e_2)=\frac{b}{2}.  $$ Therefore we obtain from
\eqref{rec652}
\begin{align}\label{rec752} \int \bar{h}(u)e^{\gamma
  C_2(y-z)}\,N_{1,2}(d(y,u,z)) =\frac{2}{\gamma^2}e^{\gamma
  v_2}-\frac{2}{\gamma^2}-\frac{2v_2}{\gamma}.  \end{align}

To evaluate $$ v_{1,1}=\int h(K,u)\,\Psi_1(K;du) $$ we split the
integration according to $u\in\{-e_1,e_1,-e_2,e_2\}$.  As all four
integrals yield the same value $ab/4$, we get $v_{1,1}=v_2$.
Summarizing, we obtain from \eqref{45} \begin{align*}
  \sigma(V_0,V_2)=\; &p(1-p)-4(1-p)^2\gamma v_2(1-1\gamma v_{2})H(\gamma
  v_2)\\ &-4(1-p)^2\big((1-p)^{-1}-1-\gamma v_2\big), \end{align*}
that is
\begin{empheq}[box=\fbox]{align}
  \label{4522}
  \sigma(V_0,V_2)=(1-p)\left[2(1-p)\gamma v_2-3p-(1-p)\gamma
  v_2(4-2\gamma v_{2})H(\gamma v_2)\right].
\end{empheq}
Figure~\ref{cov02colorall} visualizes this asymptotic covariance.

Next we turn to $\sigma(V_1,V_1)$ as given by \eqref{sigma11} for $n=2$.
Some of our calculations will also be required to compute
$\sigma(V_0,V_1)$ and $\sigma(V_0,V_0)$.  We have \begin{align}
  C_1(x;\bar{h})=\; &\frac{a}{4}\int\1\{y-x\in K^0,y\in A^+_1\cup
  A^-_1\}\,\cH^1(dy)\\ &+\frac{b}{4}\int\1\{y-x\in K^0,y\in A^+_2\cup
  A^-_2\}\,\cH^1(dy) \end{align} and a straightforward calculation (left
to the reader) yields \begin{align}\label{4525}
  C_1(x;\bar{h})=\1\{|x_1|\le a,|x_2|\le
  b\}\Big(\frac{a}{4}(b-|x_2|)+\frac{b}{4}(a-|x_1|)\Big) \end{align} as
well as \begin{align}\label{4529} C_1(x)=\frac{1}{2}\1\{|x_1|\le
  a,|x_2|\le b\}((a-|x_1|)+(b-|x_2|)).  \end{align}
From $C_1(x;\bar{h})=C_1(-x;\bar{h})$ (see \eqref{4525}) and
\eqref{rect52} (with $f(x):=e^{\gamma C_2(x)}$ and $\tilde{f}\equiv 1$)
it follows that \begin{align*} \int e^{\gamma
  C_2(y-z)}C_1(y-z;\bar{h})\,N_{1,2}(d(z,u,y))=J_1+J_2, \end{align*}
where
\begin{align*}
  J_1:=\; &\frac{a}{8}\int_{[-a,a]\times[-b,b]} e^{\gamma(a-|z_1|)(b-|z_2|)}((a-|z_1|)+(b-|z_2|))
  (b-|z_2|)\,d(z_1,z_2)
\end{align*}
and $J_2$ is defined similarly. We
have \begin{align*} J_1&=\frac{a}{2}\int \1\{0\le z_1\le a,0\le z_2\le
  b\}e^{\gamma z_1z_2}(z_1+z_2)z_2\,d(z_1,z_2)\\
  &=\frac{a^2b^2}{2}\int^1_0\int^1_0 e^{\gamma ab st}(as+bt)t\,ds\,dt\\
  &=\frac{a^3b^2}{2}\int^1_0\int^1_0 e^{\gamma ab st}st\,ds\,dt
  +\frac{a^2b^3}{2}\int^1_0\int^1_0 e^{\gamma ab st}t^2\,ds\,dt.
\end{align*} Together with the analogous formula for $J_2$ this yields
\begin{align*} \int e^{\gamma
  C_2(y-z)}C_1(y-z;\bar{h})\,N_{1,2}(d(z,u,y))
  =\frac{v_1v^2_2}{2}\int^1_0\int^1_0 e^{\gamma v_2 st}(st+t^2)\,ds\,dt.
\end{align*} Now we can use  \eqref{I15} with $r=\gamma v_2$ to obtain
\begin{align}\label{rect094}\notag \int e^{\gamma
  C_2(y-z)}&C_1(y-z;\bar{h})\,N_{1,2}(d(z,u,y))\\
  &\hspace{-5mm}=\frac{v_1}{\gamma^2}e^{\gamma v_2}-\frac{v_1}{2\gamma^3v_2}e^{\gamma
  v_2} +\frac{v_1}{2\gamma^3v_2}
  -\frac{v_1}{2\gamma^2}-\frac{3v_1v_2}{4\gamma}-\frac{v_1v_2}{2\gamma}H(\gamma
  v_2).  \end{align} Similarly, \begin{align*} \int e^{\gamma
  C_2(y-z)}&C_1(y-z)\,N_{1,2}(d(z,u,y))\\ =\; &\int \1\{0\le z_1\le a,0\le
  z_2\le b\}e^{\gamma z_1 z_2}(z_1+z_2)^2\,d(z_1,z_2)\\
  =\; &ab\int^1_0\int^1_0 e^{\gamma ab st}(as+bt)^2\,ds\,dt\\
  =\; &ab(a^2+b^2)\int^1_0\int^1_0 e^{\gamma ab st}s^2\,ds\,dt
  +2a^2b^2\int^1_0\int^1_0 e^{\gamma ab st}st\,ds\,dt.  \end{align*}
It follows from \eqref{I13} and \eqref{I14} that the latter sum equals
\begin{align*} ab&(a^2+b^2)\Big(\frac{1}{\gamma^2a^2b^2}e^{\gamma
  v_2}-\frac{1}{\gamma^3a^3b^3}e^{\gamma v_2}
  +\frac{1}{\gamma^3a^3b^3}-\frac{1}{2\gamma ab}\Big)\\
  &+2a^2b^2\Big(\frac{1}{\gamma^2a^2b^2}e^{\gamma
  ab}-\frac{1}{\gamma^2a^2b^2}-\frac{1}{\gamma ab}H(\gamma ab)
  -\frac{1}{\gamma ab}\Big).
\end{align*} Therefore
\begin{align*}
  \gamma^2&\int e^{\gamma C_2(y-z)}C_1(y-z)\,N_{1,2}(d(z,u,y))\\
  &=(a^2+b^2)\Big(\frac{1}{v_2}e^{\gamma v_2}-\frac{1}{\gamma
  v^2_2}e^{\gamma v_2} +\frac{1}{\gamma v^2_2}-\frac{\gamma}{2}\Big) +2
  e^{\gamma v_2}-2-2v_2\gamma (H(\gamma v_2)+1). \end{align*}

To proceed, we need to compute the integrals \begin{align}\label{rect61}
  I^{++}_1&:=\int \1\{y\in A^+_1,z\in A^+_1\}e^{\gamma
  C_2(y-z)}\,\cH^1(dy)\,\cH^1(dz),\allowdisplaybreaks\\ \label{rect62} I^{+-}_1&:=\int
  \1\{y\in A^+_1,z\in A^-_1\}e^{\gamma
  C_2(y-z)}\,\cH^1(dy)\,\cH^1(dz),\\ \label{rect63} I^+_{1,2}&:=\int
  \1\{y\in A^+_1,z\in A^+_2\}e^{\gamma C_2(y-z)}\,\cH^1(dy)\,\cH^1(dz)
\end{align} as well as  $I^{++}_2$ (resp.\ $I^{+-}_2$), arising from
$I^{++}_1$ (resp.\ $I^{+-}_1$) by replacing $(A^+_1,A^+_1)$ (resp.\
$(A^+_1,A^-_1)$) with $(A^+_2,A^+_2)$ (resp.\ $(A^+_2,A^-_2)$).  A
straightforward calculation gives \begin{align}\label{rect91} I^{++}_1&=
  \frac{2b}{\gamma a}e^{\gamma ab}-\frac{2}{\gamma^2a^2}e^{\gamma
  ab}+\frac{2}{\gamma^2a^2},& I^{++}_2&=\frac{2a}{\gamma b}e^{\gamma
  ab}-\frac{2}{\gamma^2b^2}e^{\gamma ab}+\frac{2}{\gamma^2b^2},\\
  \label{rect92} I^{+-}_1&=b^2,& I^{+-}_2&=a^2,\\ \label{rect93}
  I^+_{1,2}&=ab(H(\gamma ab)+1).  \end{align} We prove here
\eqref{rect91}. The proof of \eqref{rect92} and \eqref{rect93} is even
simpler.  By the parametrisation $y=(a/2,s)$ with $s \in [-b/2,b/2]$ for
$y\in A^+_1$ and $z=(a/2,t)$ with $t \in [-b/2,b/2]$ for $z\in A^+_2$ we
get \begin{align*} I^{++}_1=\int^{b/2}_{-b/2}\int^{b/2}_{-b/2} e^{\gamma
  a(b-|s-t|)}\,ds\,dt.  \end{align*} Splitting the domain of integration
into $s<t$ and $s\ge t$ yields \begin{align*} I^{++}_1&=2 e^{\gamma
  ab}\int^{b/2}_{-b/2} e^{\gamma at}\int^{b/2}_t e^{-\gamma a
  s}\,ds\,dt\\
  &=\frac{2}{\gamma a} e^{\gamma ab}\int^{b/2}_{-b/2}
  \big(1-e^{-\gamma ab/2}e^{\gamma at}\big)\,dt.  \end{align*}
Continuing this calculation gives \begin{align*}
  I^{++}_1&=\frac{2b}{\gamma a}e^{\gamma
  ab}-\frac{2}{\gamma^2a^2}e^{\gamma ab/2} \big(e^{\gamma
  ab/2}-e^{-\gamma ab/2}\big) \end{align*} and hence the first identity
in \eqref{rect91}. The second follows by symmetry.

By symmetry arguments we have \begin{align*} \int  e^{\gamma
  C_2(x-y)}\,N_{1,1}(d(x,u,y,v))=2\frac{1}{4}I^{++}_1+2\frac{1}{4}I^{++}_2
  +2\frac{1}{4}I^{+-}_1+2\frac{1}{4}I^{+-}_2+8\frac{1}{4}I^+_{1,2},
\end{align*} so that \eqref{rect91}--\eqref{rect93} yield \begin{align*}
  \int  e^{\gamma C_2(x-y)}&\,N_{1,1}(d(x,u,y,v)) =\frac{b}{\gamma
  a}e^{\gamma ab}-\frac{1}{\gamma^2a^2}e^{\gamma
  ab}+\frac{1}{\gamma^2a^2}\\ &+\frac{a}{\gamma b}e^{\gamma
  ab}-\frac{1}{\gamma^2b^2}e^{\gamma ab}+\frac{1}{\gamma^2b^2}
  +\frac{a^2+b^2}{2}+2ab(H(\gamma ab)+1).  \end{align*}
Therefore, \begin{align*} \gamma \int  &e^{\gamma
  C_2(x-y)}\,N_{1,1}(d(x,u,y,v))\\ =\; &\frac{a^2+b^2}{v_2}e^{\gamma
  v_2}-\frac{a^2+b^2}{\gamma v^2_2}e^{\gamma v_2} +\frac{a^2+b^2}{\gamma
  v^2_2} +\gamma\frac{a^2+b^2}{2}+2\gamma v_2 (H(\gamma v_2)+1),
\end{align*} so that \begin{align*} \gamma^2\int &e^{\gamma
  C_2(y-z)}C_1(y-z)\,N_{1,2}(d(z,u,y)) +\gamma \int  e^{\gamma
  C_2(x-y)}\,N_{1,1}(d(x,u,y,v))\\
  =\; &(a^2+b^2)\Big(\frac{1}{v_2}e^{\gamma v_2}-\frac{1}{\gamma
  v^2_2}e^{\gamma v_2} +\frac{1}{\gamma v^2_2}-\frac{\gamma}{2}\Big)+2
  e^{\gamma v_2}-2\\ &+\frac{a^2+b^2}{v_2}e^{\gamma
  v_2}-\frac{a^2+b^2}{\gamma v^2_2}e^{\gamma v_2} +\frac{a^2+b^2}{\gamma
  v^2_2}+\gamma\frac{a^2+b^2}{2}\\
  =\; &2(a^2+b^2)\Big(\frac{1}{v_2}e^{\gamma v_2}-\frac{1}{\gamma
  v^2_2}e^{\gamma v_2} +\frac{1}{\gamma v^2_2}\Big)+2 e^{\gamma v_2}-2.
\end{align*}
Now we can conclude from \eqref{sigma11}, \eqref{rect1} and
\eqref{rect71} that
\begin{align*}
  \sigma (V_1,V_1) =\; & 4(1-p)^2\gamma^2 v^2_1v_2H(\gamma v_2)   -4(1-p)^2v^2_1 \gamma^2 \bigg(\frac{1}{\gamma^2v_2}e^{\gamma v_2}-\frac{1}{\gamma^2v_2}-\frac{1}{\gamma}\bigg)\\
  \notag &  + (1-p)^2 2(a^2+b^2)\Big(\frac{1}{v_2}e^{\gamma v_2}-\frac{1}{\gamma v^2_2}e^{\gamma v_2} +\frac{1}{\gamma v^2_2}\Big)+(1-p)^2(2 e^{\gamma v_2}-2),
\end{align*} that is
\begin{empheq}[box=\fbox]{align}
  \label{4527}
  \begin{aligned}
  \sigma (V_1,V_1) =\; &(1-p) \bigg[ 2p+4(1-p)\gamma^2 v^2_1v_2H(\gamma
    v_2) \bigg. \\  &\bigg.-4\gamma^2v^2_1
    \bigg(\frac{p}{\gamma^2v_2}-\frac{1-p}{\gamma}\bigg) + 2
  (a^2+b^2)\Big(\frac{1}{v_2}-\frac{p}{\gamma v^2_2}\Big)\bigg],
  \end{aligned}
\end{empheq}
which is shown in Fig.~\ref{var1colorall}.

We now turn to $\sigma(V_0,V_1)$.  It follows from \eqref{rect094} that
\begin{align}\label{rect68} -2\gamma^3&\int e^{\gamma
  C_2(y-z)}C_1(y-z;\bar{h})\,N_{1,2}(d(z,u,y))\\ \notag =\; &-2\gamma
  v_1e^{\gamma v_2}+\frac{v_1}{v_2}e^{\gamma v_2} -\frac{v_1}{v_2}
  +\gamma v_1+\frac{3\gamma^2v_1v_2}{2}+\gamma^2 v_1v_2H(\gamma v_2).
\end{align}

Furthermore we have from \eqref{rect71} and $v_{1,1}=v_2$
\begin{align*}
  (\gamma^3v_{1,1}&-\gamma^2)\int e^{\gamma
  C_2(y-z)}\,N_{1,2}(d(z,u,y))\\
  &=2v_1\gamma^3v_2\bigg(\frac{1}{\gamma^2v_2}e^{\gamma
  v_2}-\frac{1}{\gamma^2v_2}-\frac{1}{\gamma}\bigg)
  -2v_1\gamma^2\bigg(\frac{1}{\gamma^2v_2}e^{\gamma
  v_2}-\frac{1}{\gamma^2v_2}-\frac{1}{\gamma}\bigg)\allowdisplaybreaks\\ &=2\gamma
  v_1e^{\gamma v_2}-2\gamma v_1-2\gamma^2v_1v_2
  -\frac{2v_1}{v_2}e^{\gamma v_2}+\frac{2v_1}{v_2}+2\gamma v_1.
\end{align*} Further we have from \eqref{rec752} and $\overline\Psi(1)=v_1$
\begin{align*} 2\gamma^3\overline{\Psi}_1(1)&\int e^{\gamma
  C_2(y-z)}\bar{h}(u)\,N_{1,2}(d(z,u,y))
  =2\gamma^3v_1\Big(\frac{2}{\gamma^2}e^{\gamma
  v_2}-\frac{2}{\gamma^2}-\frac{2v_2}{\gamma}\Big)\\ &=4\gamma
  v_1e^{\gamma v_2}-4\gamma v_1-4\gamma^2v_1v_2.  \end{align*}
Summarizing the previous formulae we arrive at
\begin{align}\label{rect85}\notag
  (1-p)^2&\int\big(\chi'(y-z)+2\gamma^3\overline{\Psi}_1(1) e^{\gamma
  C_2(y-z)}\bar{h}(u)\big)\,N_{1,2}(d(z,u,y))\\ 
  =\; & (1-p)\gamma v_1\left[ 1 + \left(3-\frac{1}{\gamma v_2}\right)p + (1-p)\gamma v_2 \left( H(\gamma v_2) -\frac{9}{2} \right)\right].
\end{align}

Next we consider \begin{align*} I:=\int e^{\gamma
  C_2(y-z)}\bar{h}(u)\,N_{1,1}(d(y,u,z,v)).  \end{align*} Then
$I=I_1+I_2$, where \begin{align*} I_1&:=\frac{a}{8}\int \1\{y\in
  A^+_1\cup A^-_1\}e^{\gamma C_2(y-z)}\,\cH^1(dy)\,\cH^1(dz),\\
  I_2&:=\frac{b}{8}\int \1\{y\in A^+_2\cup A^-_2\}e^{\gamma
  C_2(y-z)}\,\cH^1(dy)\,\cH^1(dz).  \end{align*} By symmetry, $$
I_1=2\frac{a}{8}I^{++}_1+2\frac{a}{8}I^{+-}_1+4\frac{a}{8}I^+_{1,2},
\quad
I_2=2\frac{b}{8}I^{++}_2+2\frac{b}{8}I^{+-}_2+4\frac{a}{8}I^+_{1,2}, $$
where the occurring integrals have been defined by
\eqref{rect61}--\eqref{rect63}.  The formulae
\eqref{rect91}--\eqref{rect93} yield \begin{align*}
  I=\; &\Big(\frac{a}{2\gamma}+\frac{b}{2\gamma}\Big)e^{\gamma ab}
  -\Big(\frac{1}{2\gamma^2a}+\frac{1}{2\gamma^2b}\Big)e^{\gamma
  ab}+\frac{1}{2\gamma^2a}+\frac{1}{2\gamma^2b}\\
  &+\frac{ab^2}{4}+\frac{a^2b}{4}+\Big(\frac{a^2b}{2}+\frac{ab^2}{2}\Big)(H(\gamma
  ab)+1), \end{align*} that is \begin{align}\notag 
  I=\frac{v_1}{2\gamma}e^{\gamma v_2}-\frac{v_1}{2\gamma^2v_2}e^{\gamma
  v_2}
  +\frac{v_1}{2\gamma^2v_2}+\frac{3v_1v_2}{4}+\frac{v_1v_2}{2}H(\gamma
  v_2).  \end{align}
It follows that
\begin{align}\label{rect595}\notag
  -2(1-p)^2\gamma^2&\int e^{\gamma
  C_2(y-z)}\bar{h}(u)\,N_{1,1}(d(y,u,z,v))\\
  =\; & (1-p)\gamma v_1\left[ \frac{p}{\gamma v_2} -1 -\left( \frac{3}{2} +H(\gamma
  v_2)\right)(1-p)\gamma v_2 \right].
\end{align}
Now we conclude
from \eqref{46} and \eqref{rect1} that
\begin{align}
  \sigma(V_0,V_1)&=
  (1-p)\gamma v_1\left[ 1-p + 4(1-p)\gamma v_2(1-\gamma v_2) H(\gamma v_2) \right]+c_{1,2}+c_{1,1},
\end{align}
where $c_{1,2}$ is given by the right-hand side of
\eqref{rect85} and $c_{1,1}$ is given by the right-hand side of
\eqref{rect595}. Thus, we derive
\begin{empheq}[box=\fbox]{align}
  \label{sigma01}
  \sigma(V_0,V_1) = (1-p)\gamma v_1\left[ 1+2p  + (1-p)\gamma v_2\left(4(1-\gamma v_2)H(\gamma v_2) -6\right) \right].
\end{empheq}
The asymptotic covariance is plotted in Fig.~\ref{cov01colorall}.

Finally, we determine $\sigma(V_0,V_0)$, as given by \eqref{47}.
From \eqref{rect52} and \eqref{4525} we get \begin{align*} \int& e^{\gamma
  C_2(y-z)}C_1(y-z;\bar{h})\bar{h}(u) \, N_{1,2}(d(z,u,y))\\
  =\; &\frac{a}{4}\int \1\{|z_1|\le a,|z_2|\le b\}e^{\gamma
  |z_1||z_2|}|z_2|
  \Big(\frac{a}{4}|z_2|+\frac{b}{4}|z_1|\Big)\,d(z_1,z_2)\\
  &+\frac{b}{4}\int \1\{|z_1|\le a,|z_2|\le b\}e^{\gamma
  |z_1||z_2|}|z_1|
  \Big(\frac{a}{4}|z_2|+\frac{b}{4}|z_1|\Big)\,d(z_1,z_2)\\
  =\; &\frac{a}{4}\int^b_0\int^a_0 e^{\gamma z_1
  z_2}z_2(az_2+bz_1)\,dz_1\,dz_2 +\frac{b}{4}\int^b_0\int^a_0e^{\gamma
  z_1z_2}z_1 (az_2+bz_1)\,dz_1\,dz_2\\
  =\; &\frac{a^2b^2}{4}\int^1_0\int^1_0 e^{\gamma abst}t(abt+abs)\,ds\,dt
  +\frac{a^2b^2}{4}\int^1_0\int^1_0e^{\gamma abst}s (abt+abs)\,ds\,dt\\
  =\; &\frac{a^3b^3}{2}\int^1_0\int^1_0e^{\gamma abst}(s^2+st)\,ds\,dt.
\end{align*} Using now \eqref{I15}, we obtain
\begin{align*} 4\gamma^4\int& e^{\gamma
  C_2(y-z)}C_1(y-z;\bar{h})\bar{h}(u) \, N_{1,2}(d(z,u,y))\\
  &=4\gamma^2v_2e^{\gamma v_2}-2\gamma e^{\gamma v_2}+2\gamma
  -2\gamma^2v_2-3\gamma^3v^2_2-2\gamma^3v^2_2 H(\gamma v_2).
\end{align*} From \eqref{rec752} \begin{align*}
  (-4\gamma^4v_2+4\gamma^3)\int& e^{\gamma C_2(y-z)}\bar{h}(u) \,
  N_{1,2}(d(z,u,y))\\
  &=(-4\gamma^4v_2+4\gamma^3)\Big(\frac{2}{\gamma^2}e^{\gamma
  v_2}-\frac{2}{\gamma^2}-\frac{2v_2}{\gamma}\Big)\\ &=-8\gamma^2v_2
  e^{\gamma v_2}+8\gamma^2v_2+8\gamma^3v^2_2+8\gamma e^{\gamma
  v_2}-8\gamma-8\gamma^2v_2.  \end{align*}
Summarizing, we have
\begin{align}\label{d12} \notag
  (1-p)^2\int& \bar{h}(u)\chi''(y-z) \, N_{1,2}(d(z,u,y))\\
  =\; & (1-p) \gamma \left[ 6p - 4\gamma v_2 + (1-p) \gamma v_2(5\gamma v_2 -
  2H(\gamma v_2)\gamma v_2 -2) \right]
\end{align}

Next we note that \begin{align*} \int e^{\gamma
  C_2(y-z)}&\bar{h}(u)\bar{h}(v)\,N_{1,1}(d(y,u,z,v))\\
  &=2\frac{a^2}{16} I^{++}_1+2\frac{a^2}{16}
  I^{+-}_1+2\frac{b^2}{16}I^{++}_2
  +2\frac{b^2}{16}I^{+-}_2+8\frac{ab}{16}I^+_{1,2}, \end{align*} where
$I^{++}_1,I^{+-}_1,I^+_{1,2},I^{++}_2,I^{+-}_2$ have been defined by
\eqref{rect61}--\eqref{rect63}.  The formulae
\eqref{rect91}--\eqref{rect93} give \begin{align*} \int & e^{\gamma
  C_2(y-z)}\bar{h}(u)\bar{h}(v)\,N_{1,1}(d(y,u,z,v))\\
  =\; &\frac{a^2}{8}\Big(\frac{2b}{\gamma a}e^{\gamma
  ab}-\frac{2}{\gamma^2a^2}e^{\gamma ab}+\frac{2}{\gamma^2a^2}\Big)
  +\frac{a^2b^2}{8} +\frac{b^2}{8}\Big(\frac{2a}{\gamma b}e^{\gamma
  ab}-\frac{2}{\gamma^2b^2}e^{\gamma ab}+\frac{2}{\gamma^2b^2}\Big)\\
  &+\frac{a^2b^2}{8}+\frac{a^2b^2}{2}(H(\gamma ab)+1)\\
  =\; &\frac{v_2}{2\gamma}e^{\gamma v_2}-\frac{1}{2\gamma^2}e^{\gamma
  v_2}+\frac{1}{2\gamma^2} +\frac{v^2_2}{4}
+\frac{v^2_2}{2}H(\gamma v_2)+\frac{v^2_2}{2}.  \end{align*}
Therefore
\begin{align}\label{d11}\notag 4(1-p)^2\gamma^3&\int  e^{\gamma
  C_2(y-z)}\bar{h}(u)\bar{h}(v)\,N_{1,1}(d(y,u,z,v))\\
  =\; & 
  (1-p)\gamma\left[ (2H(\gamma v_2) + 3)(1-p)(\gamma v_2)^2 + 2\gamma v_2 -2p \right].
\end{align}
Now we conclude from \eqref{47} and \eqref{rect1} that
\begin{align}
  \sigma(V_0,V_0)=(1-p)\gamma\left[\right. &1-2p+(2p-3)\gamma v_2\\
    &+\left. 4(1-p)(1-\gamma v_2)^2\gamma v_2 H(\gamma v_2) \right] +d_{1,2}+d_{1,1},
\end{align}
where $d_{1,2}$ is given by the right-hand side of \eqref{d12} and
$d_{1,1}$ is given by the right-hand side of \eqref{d11}. Thus, we finally
derive 
\begin{empheq}[box=\fbox]{align}
  \label{rect6789}
  \begin{aligned}
   \sigma(V_0,V_0)=(1-p)\gamma\left[\right. &1+2p+(4p-7)\gamma v_2\\
    &+\left. 4(1-p)\gamma v_2\left(2\gamma v_2 +(1-\gamma v_2)^2
    H(\gamma v_2)\right)
  \right],   
\end{aligned}
\end{empheq}
which is plotted in Fig.~\ref{var0colorall}.

The reader might have noticed that the asymptotic covariances
$\sigma(V_i,V_j)$ (with
the exception of $\sigma(V_1,V_1)$) depend on the parameters $\gamma$, $v_1$, and $v_2$
in a specific way.
In order to explain these invariance properties, let $Z_{a,b,\gamma}$
denote the Boolean model with grains $K=[0,a]\times[0,b]$ and intensity $\gamma$.
By applying to each rectangle
of the underlying Poisson process the linear transformation
$T_{a,b}: \R^2\to\R^2 ,(x_1,x_2)\mapsto (ax_1,bx_2)$,
one obtains the distributional identity
$$
Z_{a,b,\gamma} \overset{d}{=} T_{a,b} Z_{1,1,ab\gamma}.
$$
Together with the fact that $V_i(T_{a,b}A)=(ab)^{i/2} V_i(A)$, for $A\in\mathcal{R}^2$ and $i\in\{0,2\}$, we see that
\begin{align*}
\sigma(V_i,V_j) & =\lim_{r(W)\to\infty}\frac{\cov(V_i(Z_{a,b,\gamma}\cap W),V_j(Z_{a,b,\gamma}\cap W))}{V_2(W)}\\
 & =\lim_{r(W)\to\infty}\frac{\cov(V_i(T_{a,b}Z_{1,1,ab\gamma}\cap W),V_j(T_{a,b}Z_{1,1,ab\gamma}\cap W))}{V_2(W)}\\
 & =(ab)^{i/2+j/2-1} \lim_{r(W)\to\infty}\frac{\cov(V_i(Z_{1,1,ab\gamma}\cap T_{a,b}^{-1} W),V_j(Z_{1,1,ab\gamma}\cap T_{a,b}^{-1}W))}{V_2(T_{a,b}^{-1}W)}\\
  & =v_2^{i/2+j/2-1} \lim_{r(W)\to\infty}\frac{\cov(V_i(Z_{1,1,\gamma v_2}\cap W),V_j(Z_{1,1,\gamma v_2}\cap W))}{V_2(W)},
\end{align*}
for $i,j\in\{0,2\}$.
This shows that for all Boolean models of deterministic rectangles with
fixed $\gamma v_2$, the asymptotic covariances between volume and Euler
characteristic are a power of $v_2$ times a constant depending on
$\gamma v_2$.

Next we investigate the invariance properties of $\sigma(V_0,V_1)$ and $\sigma(V_1,V_2)$. For $i\in\{1,2\}$ we define
$$
V_{1,e_i}(A):=\int {\bf 1}\{u=\pm e_i\} \, \Psi_1(A;du), \quad A\in\mathcal{R}^2,
$$
which are again geometric functionals. If $W$ is a rectangle with sides in the directions $e_1$ and $e_2$, which we can assume in the following, we have that
$$
V_1(Z_{a,b,\gamma}\cap W)=V_{1,e_1}(Z_{a,b,\gamma}\cap W)+V_{1,e_2}(Z_{a,b,\gamma}\cap W).
$$
By the same arguments as in the previous computation we obtain that, for $i\in\{0,2\}$,
\begin{align*}
\sigma(V_i,V_1) & =\lim_{r(W)\to\infty}\frac{\cov(V_i(Z_{a,b,\gamma}\cap W),V_1(Z_{a,b,\gamma}\cap W))}{V_2(W)}\\
& =\lim_{r(W)\to\infty}\left\{\frac{\cov(V_i(T_{a,b}Z_{1,1,ab\gamma}\cap W),V_{1,e_1}(T_{a,b}Z_{1,1,ab\gamma}\cap W))}{V_2(W)}\right.\\
& \quad \quad \quad \quad \left.  +\frac{\cov(V_i(T_{a,b}Z_{1,1,ab\gamma}\cap W),V_{1,e_2}(T_{a,b}Z_{1,1,ab\gamma}\cap W))}{V_2(W)}\right\}\\
& = \lim_{r(W)\to\infty}\left\{(ab)^{i/2-1}a\frac{\cov(V_i(Z_{1,1,ab\gamma}\cap W),V_{1,e_1}(Z_{1,1,ab\gamma}\cap W))}{V_2(W)}\right.\\
& \quad \quad \quad \quad \left.+(ab)^{i/2-1}b\frac{\cov(V_i(Z_{1,1,ab\gamma}\cap W),V_{1,e_2}(Z_{1,1,ab\gamma}\cap W))}{V_2(W)}\right\}.
\end{align*}
Using that the asymptotic covariances between $V_i$ and $V_{1,e_1}$ and between $V_i$ and  $V_{1,e_2}$ are the same for the Boolean model $Z_{1,1,ab\gamma}$ due to symmetry, we conclude that
\begin{align*}
\sigma(V_i,V_1) & =\lim_{r(W)\to\infty}(ab)^{i/2-1}(a+b)\frac{\cov(V_i(Z_{1,1,ab\gamma}\cap W),V_{1,e_1}(Z_{1,1,ab\gamma}\cap W))}{V_2(W)}\\
& =\lim_{r(W)\to\infty}\left\{(ab)^{i/2-1}\frac{(a+b)}{2}\frac{\cov(V_i(Z_{1,1,ab\gamma}\cap W),V_{1,e_1}(Z_{1,1,ab\gamma}\cap W))}{V_2(W)}\right.\\
& \quad \quad \quad \quad \left.+(ab)^{i/2-1}\frac{(a+b)}{2}\frac{\cov(V_i(Z_{1,1,ab\gamma}\cap W),V_{1,e_2}(Z_{1,1,ab\gamma}\cap W))}{V_2(W)}\right\}\\
& = \lim_{r(W)\to\infty} v_2^{i/2-1}\frac{v_1}{2}\frac{\cov(V_i(Z_{1,1,\gamma v_2}\cap W),V_{1}(Z_{1,1,\gamma v_2}\cap W))}{V_2(W)}.
\end{align*}
Thus, the asymptotic covariance between volume and surface area is
$v_1$ times a constant depending on $\gamma v_2$, while the covariance
between Euler characteristic and surface area is $v_1 v_2^{-1}$
times a constant depending on $\gamma v_2$.

Figure~\ref{fig_cov_numerical_integration} summarizes the results of
this section visually. It shows the asymptotic covariances
$\sigma(V_i,V_j)$ as a function of the intensity $\gamma$ for Boolean
models of aligned rectangles for a variety of  aspect ratios $b/a$.
   
\section{Simulations of Boolean models with isotropic or aligned
rectangles}\label{sec:simulations}

Planar Boolean models with either squares or rectangles with aspect
ratio $1/2$ as grains are simulated in a finite observation window. We study the
variances and covariances of the intrinsic volumes as well as their
relative frequency histograms weighted by the size of each bin. We
compare the simulation results for aligned rectangles to the analytic
formulae for the covariances in the previous
Section~\ref{sec:rectangles}. Moreover, we simulate rectangles with a
uniform (isotropic) orientation distribution and find, e.g., for
$\sigma(V_0,V_1)$ a qualitatively different behaviour.

Tensor valuation densities and the density of the Euler characteristic
$\overline{V}_0$ of anisotropic Boolean models are studied
in~\cite{HoerrmannHugKlattMecke}. The same simulation procedure is
applied here with even better statistics for reliable estimates of the
second moments and the histograms of the intrinsic volumes.

The grain centres are random points, uniformly
distributed within the simulation box. The union of the rectangles is
computed using the Computational Geometry Algorithms Library
(CGAL)~\cite{CGAL}. The program \textsc{papaya} then calculates the
Minkowski functionals of the Boolean
model~\cite{SchroederTurketal2010JoM}.

For aligned rectangles, the covariances $\sigma(V_2,V_2)$,
$\sigma(V_0,V_0)$, and $\sigma(V_0,V_2)$ as well as the rescaled
covariances $\sigma(V_1,V_2)/(2a+2b)$ and $\sigma(V_0,V_1)/(2a+2b)$ are
only functions of $v_2$ and $\gamma v_2$, as shown in the previous
Section~\ref{sec:rectangles}. In other words, if the unit of area is
chosen to be the area of a single grain $v_2 = a\cdot b = 1$ (so that
the area of the typical grain does not depend on the aspect ratio), 
the rescaled covariances are independent of the aspect ratio.
Therefore, we define in the following the unit of length by the square
root of the area of a single grain.

Parts of this section are taken from the PhD thesis of one of the authors~\cite{Klatt2016}.

\subsection{Variances and covariances} \label{sec:var-cov-simulation}

The first moments of area or perimeter of a Boolean model are rather
insensitive to the grain distribution. Indeed, if the unit of area is
chosen to be the mean area of a single grain, the density of the area,
i.e., the occupied area fraction, of the Boolean model is only a
function of the intensity.
Moreover, if the density of the perimeter in the asymptotic limit is
divided by the mean perimeter of a single grain, it is also independent
of the grain distribution~\cite[Theorem 9.1.4.]{SchneiderWeil}.

Does the same hold for the second moments? Is there a qualitatively
different behaviour in the variances and covariances depending on
whether the orientation distribution of the grains is isotropic or
anisotropic? Which covariances or variances are invariant under affine
transformations of the grain distributions?

Depending on the computational costs, for each different set
of parameters we perform between $M_s=\num{21000}$ and $\num{600000}$ simulations
of Boolean models with rectangles: at varying intensities $\gamma$, with
aspect ratio $1$ or $1/2$, and for rectangles either aligned w.r.t. the
observation window or with an isotropic orientation distribution. The
simulation box is a square with side length $L=4a$, where periodic
boundary conditions are applied. The number of grains within the
simulation box is a random number and follows a Poisson distribution
with mean $\gamma\cdot L^2$.
To estimate the covariances, we simulate more than \num{5800000} samples
of Boolean models including about \num{54000000} rectangles in total.

Because of the periodic boundary conditions, the covariances of this
system coincide with the asymptotic covariances for the infinite volume
system from Section~\ref{sec:rectangles}, as we have pointed out in
Subsection~\ref{subsec:torus}.

\begin{figure}[t] \centering
  \subfigure{\includegraphics[width=0.48\textwidth]{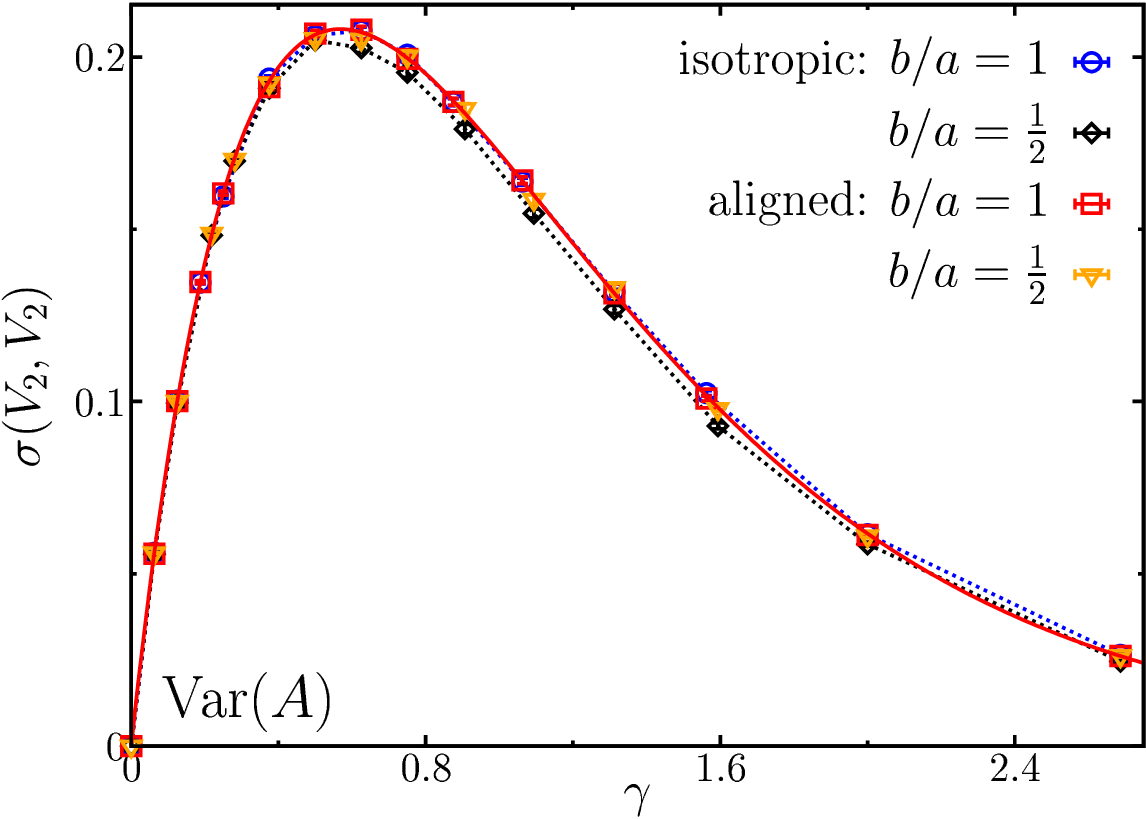}}\hfill
  \subfigure{\includegraphics[width=0.48\textwidth]{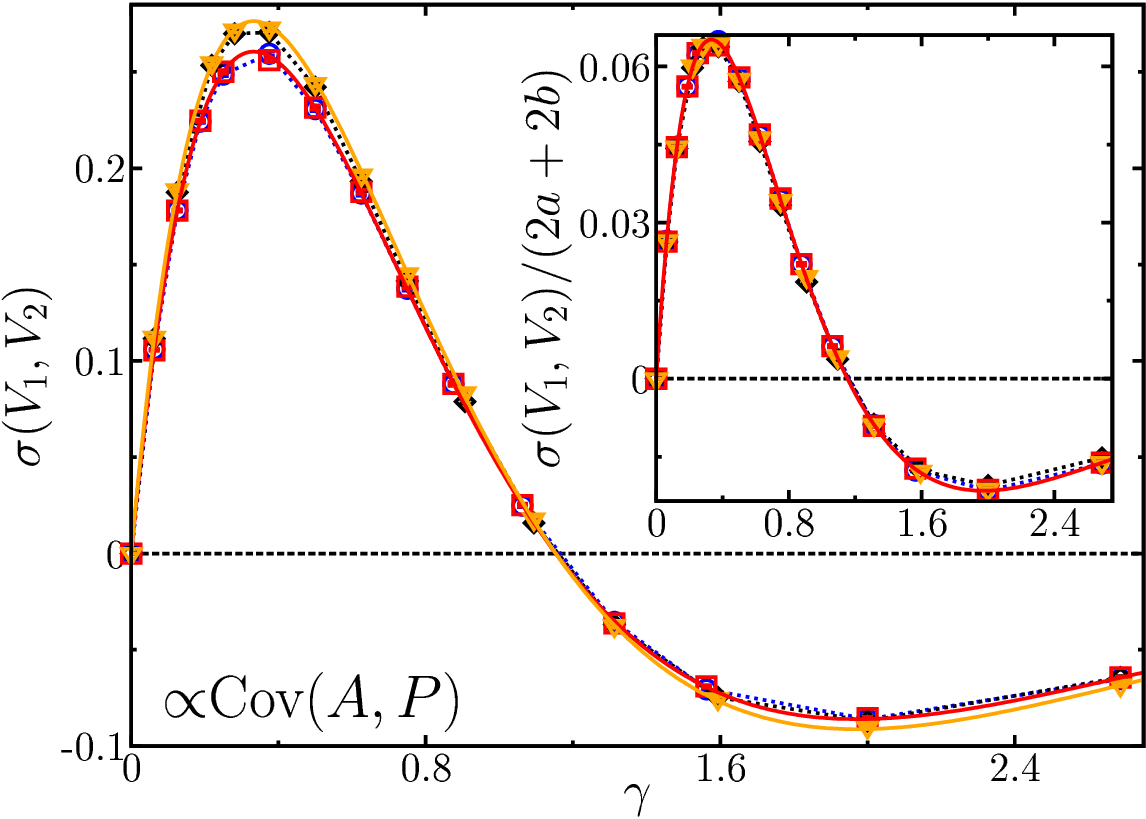}}\\
  \subfigure{\includegraphics[width=0.48\textwidth]{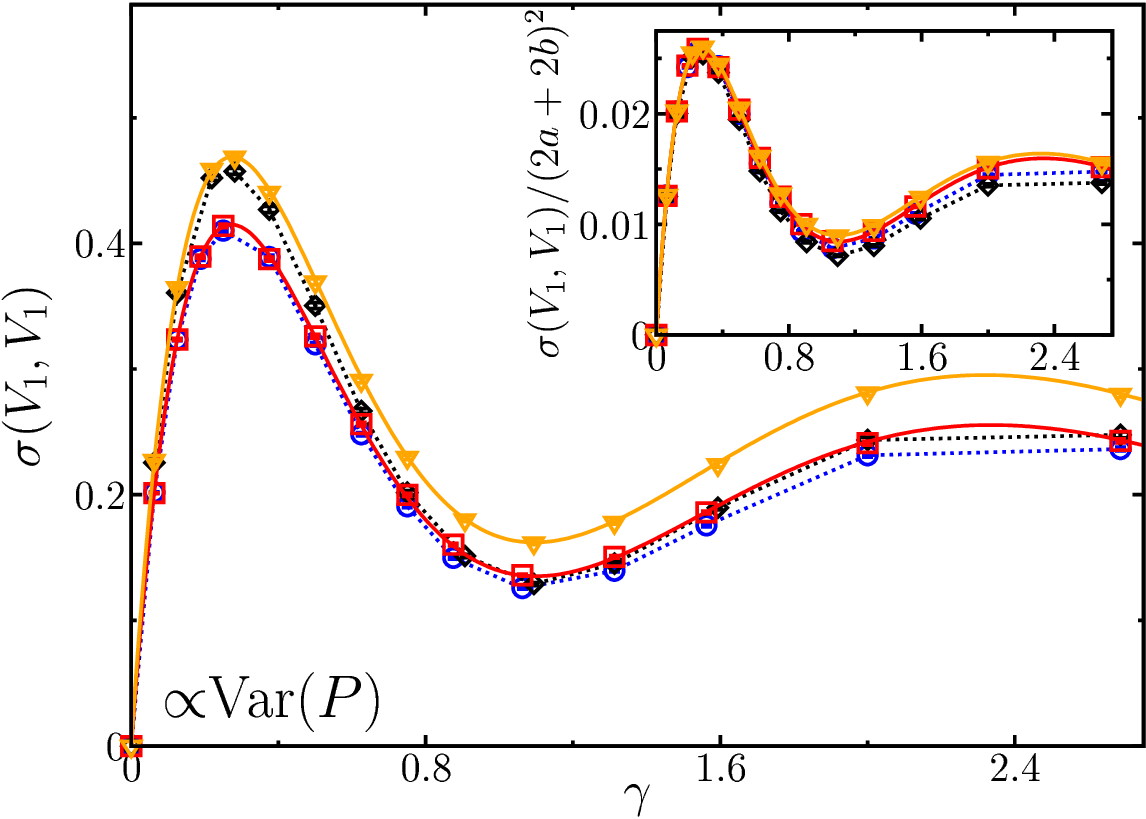} \label{fig_analytic_numeric_11}}\hfill
  \subfigure{\includegraphics[width=0.48\textwidth]{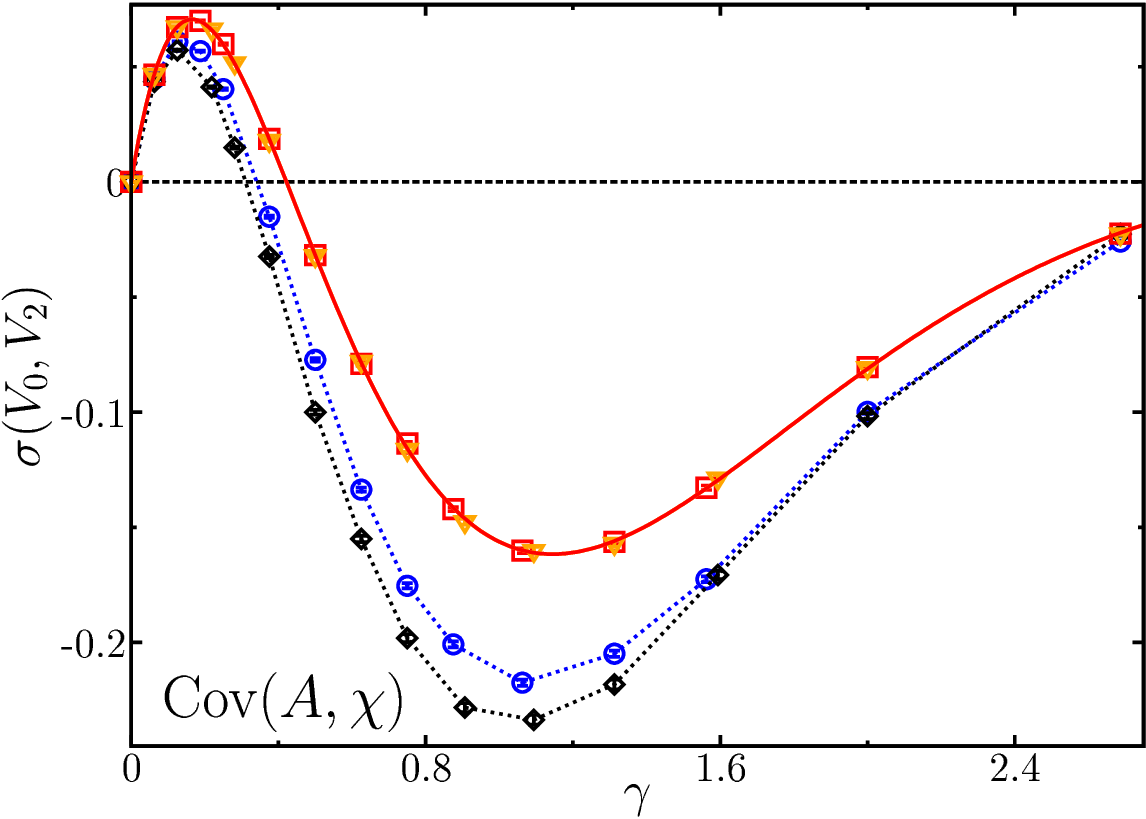}}\\
  \subfigure{\includegraphics[width=0.48\textwidth]{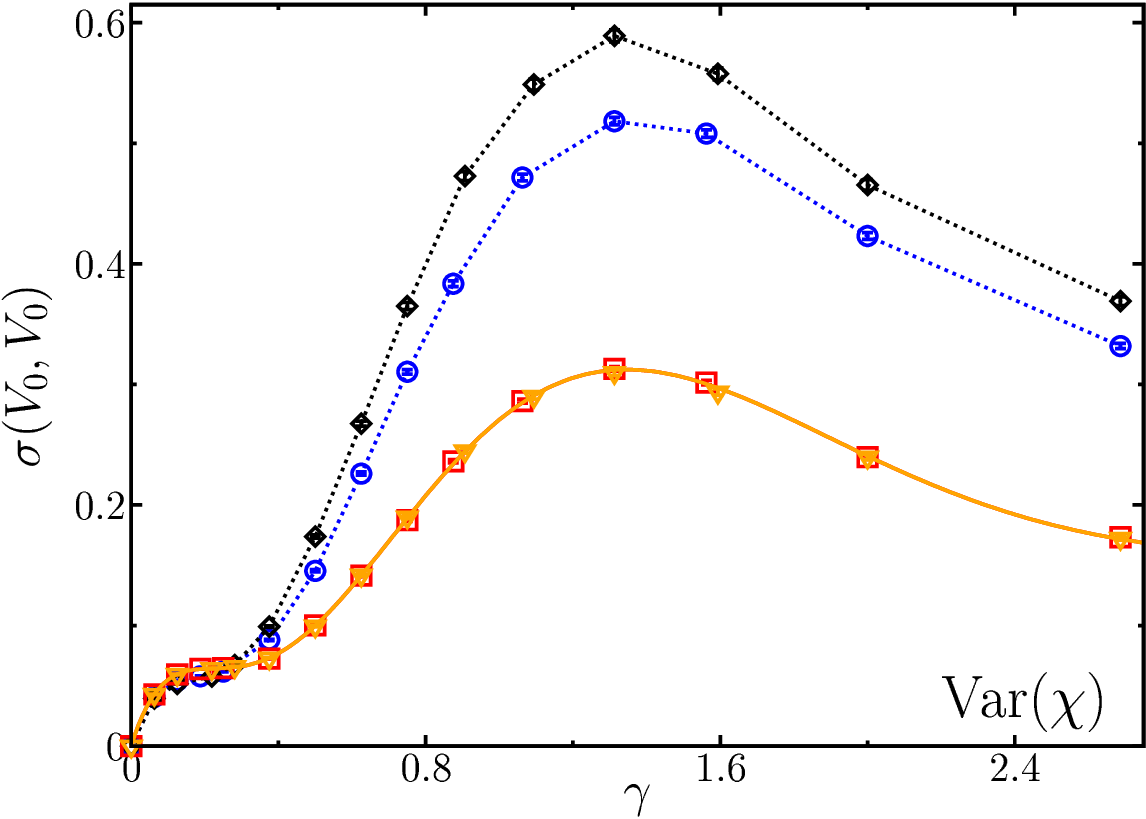}}\hfill
  \subfigure{\includegraphics[width=0.48\textwidth]{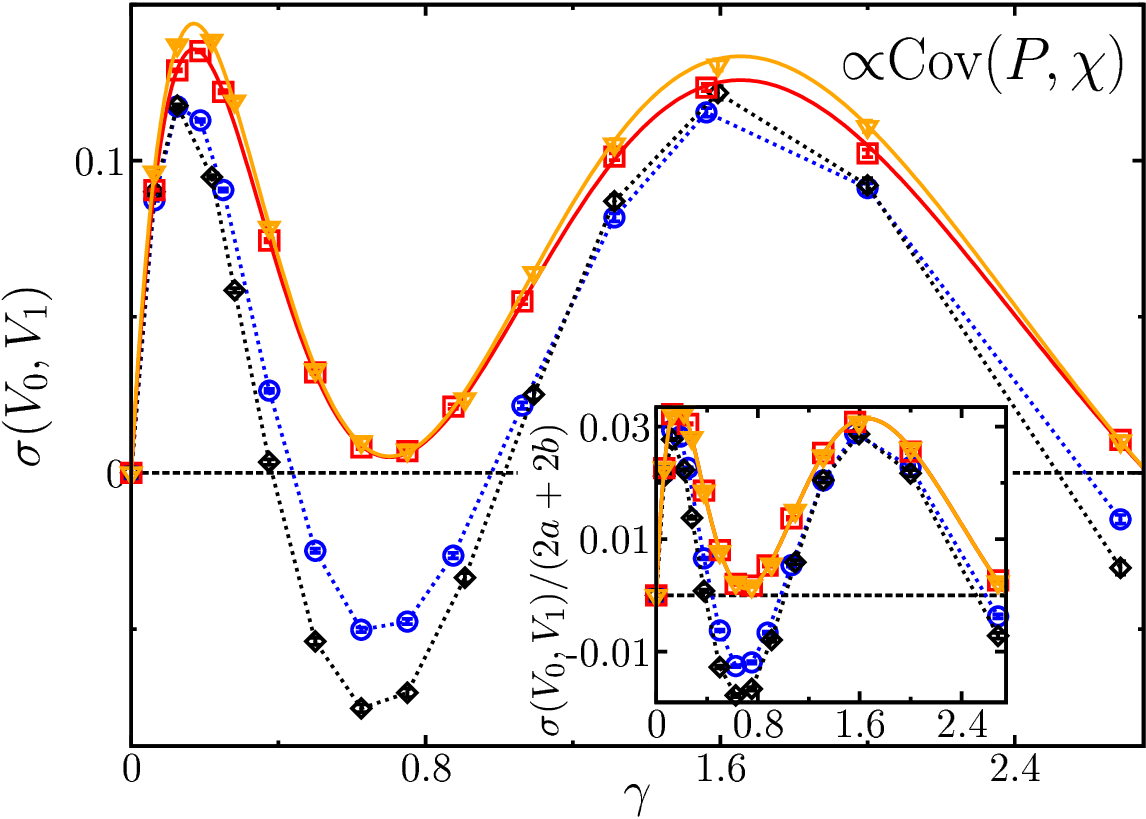}}
  \caption{Variances and covariances of the intrinsic volumes $V_2$
    (area $A$), $V_1$ (proportional to perimeter $P$), and $V_0$ (Euler
    characteristic $\chi$) of Boolean models as a function of the
    intensity $\gamma$. Depicted are both numerical estimates in finite
    observation windows with periodic boundary conditions (marks with
    dotted lines as guides to the eye) and analytic curves of the
    covariances (solid lines), see Eqs.~\protect\eqref{rect2},
    \protect\eqref{rect5}, \protect\eqref{4522}, \protect\eqref{4527}, \protect\eqref{sigma01}, and
    \protect\eqref{rect6789}.
    Four different Boolean models are simulated: both for squares
    (${b}/{a}=1$) and rectangles (${b}/{a}={1}/{2}$) either an isotropic
    orientation distribution is used or the grains are aligned with the
    $x$-direction.
    In the insets, the covariances and the variance of the perimeter of
    the Boolean model are rescaled by the perimeter of a single grain.
    In contrast to Fig.~\protect\ref{fig_cov_numerical_integration}, the unit of
  area is the size of a single grain, that is $v_2=ab=1$.} \label{fig_Boolean_cov}
\end{figure}

For each sample $m\in\{1,\dots M_s\}$ of a Boolean model, we determine
the intrinsic volumes $V_{i}^{(m)}$ ($i=0,1,2$). The sample covariance
then provides an estimate for the covariance between the Minkowski
functionals:
\[
  s({V_{i},V_{j}}) :=
  \frac{1}{M_s-1}\sum_{m=1}^{M_s}(V_{i}^{(m)}-\langle{}V_i{}\rangle)(V_{j}^{(m)}-\langle{}V_j{}\rangle)
\]
using the sample mean
$\langle{}V_i{}\rangle:=\frac{1}{M_s}\sum_{m=1}^{M_s}V_{i}^{(m)}$ as an
estimator for the expectation. In accordance with the definition of the
asymptotic covariances in Theorem~\ref{thm:CovariancesGeneral}, the
sample covariance is then divided by the size $L^2$ of the observation
window. We finally use bootstrapping (with \num{1000} bootstrap samples)
to estimate the mean and the error of the estimators.

Figure~\ref{fig_Boolean_cov} shows the simulation results for the
variances and covariances of the intrinsic volumes for an
isotropic orientation distribution of the grains as well as for aligned rectangles. In the
latter case, the simulation results are compared to the analytic results in Eqs.~\eqref{rect2},
\eqref{rect5}, \eqref{4522}, \eqref{4527}, \eqref{sigma01}, and
\eqref{rect6789}.
They are in excellent agreement.

\begin{figure}[p] \centering \hfill
  \includegraphics[width=0.45\textwidth]{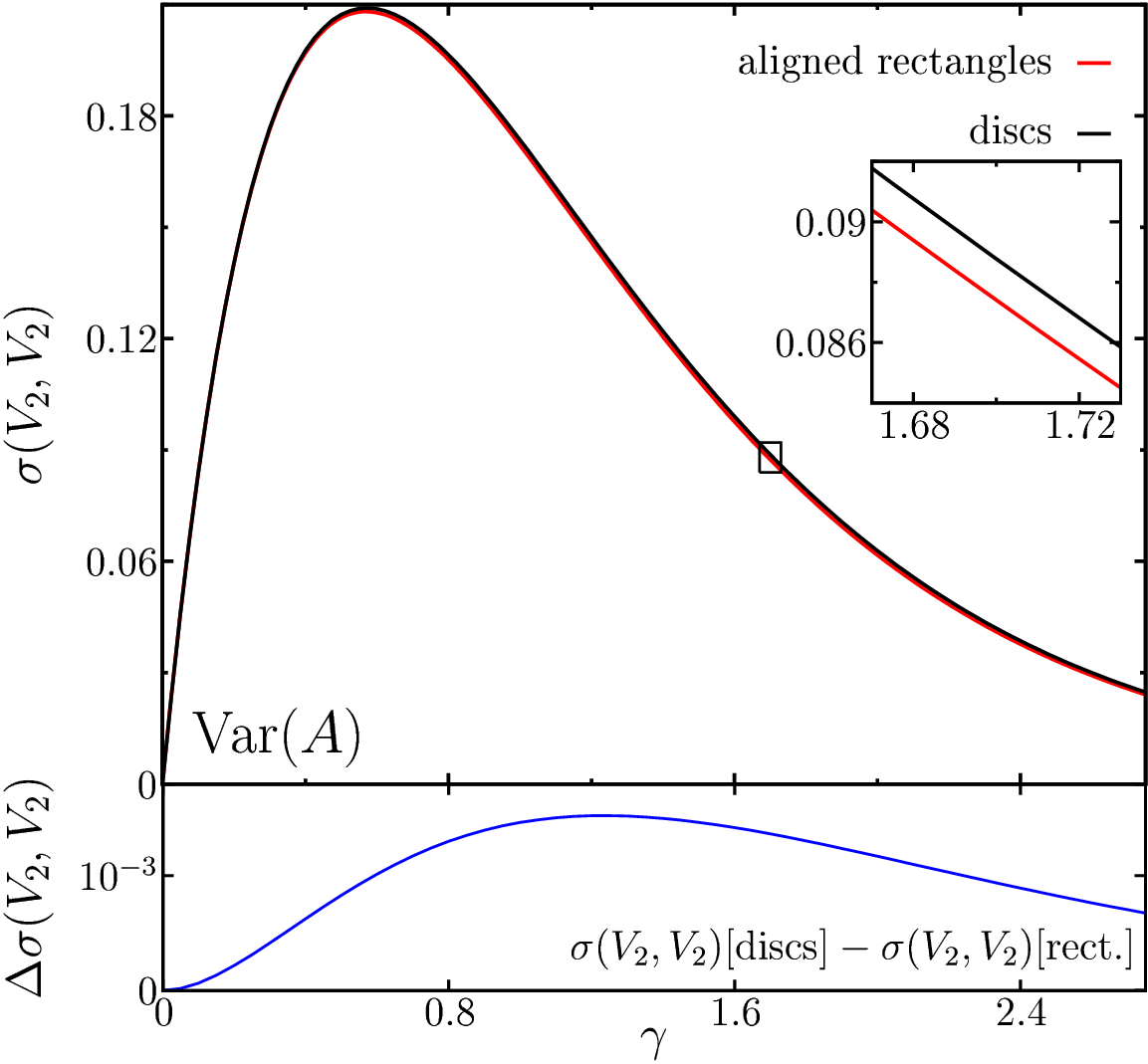}
  \hfill
  \includegraphics[width=0.45\textwidth]{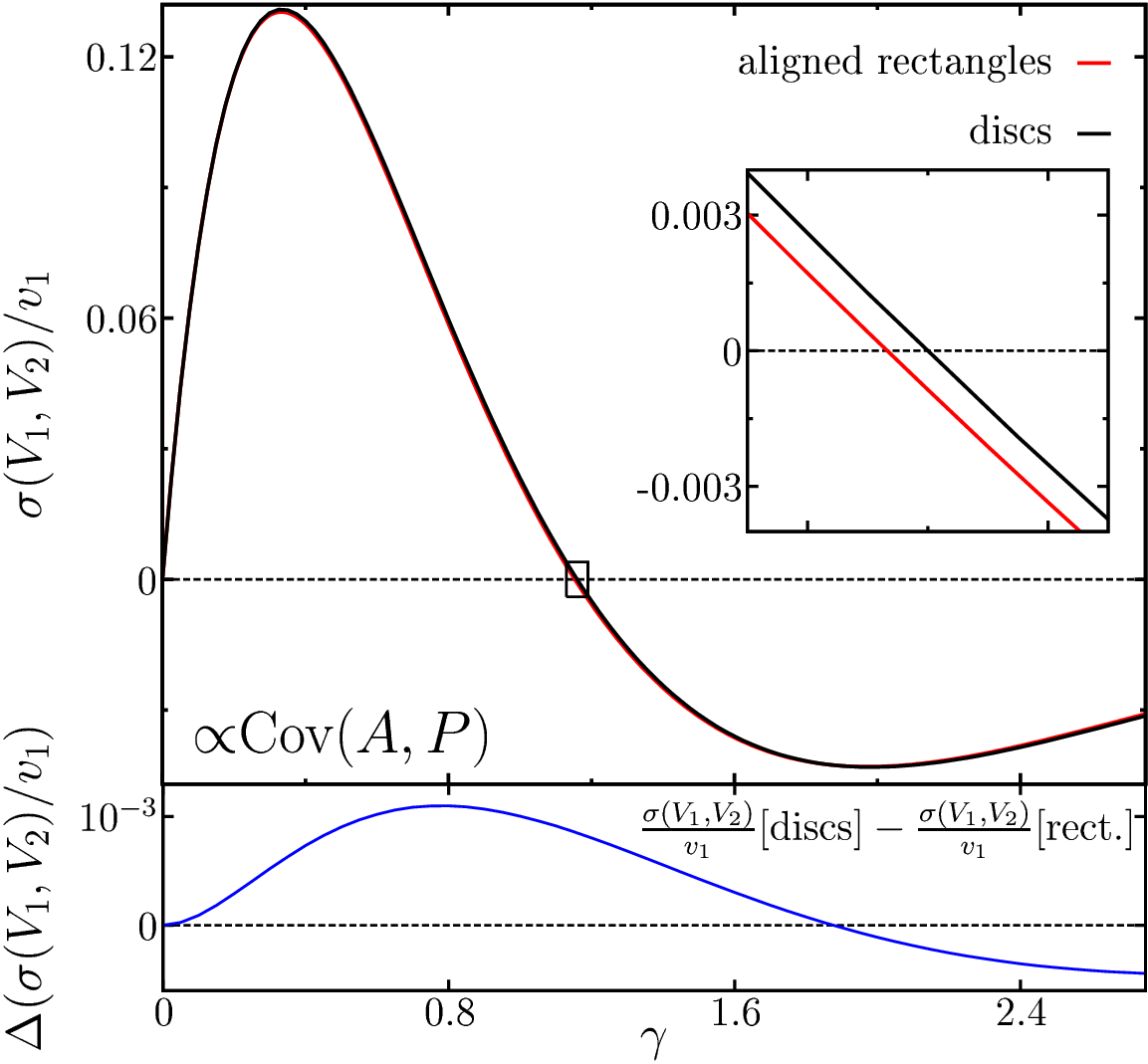}
  \hfill\hbox{}\\ \hfill
  \includegraphics[width=0.45\textwidth]{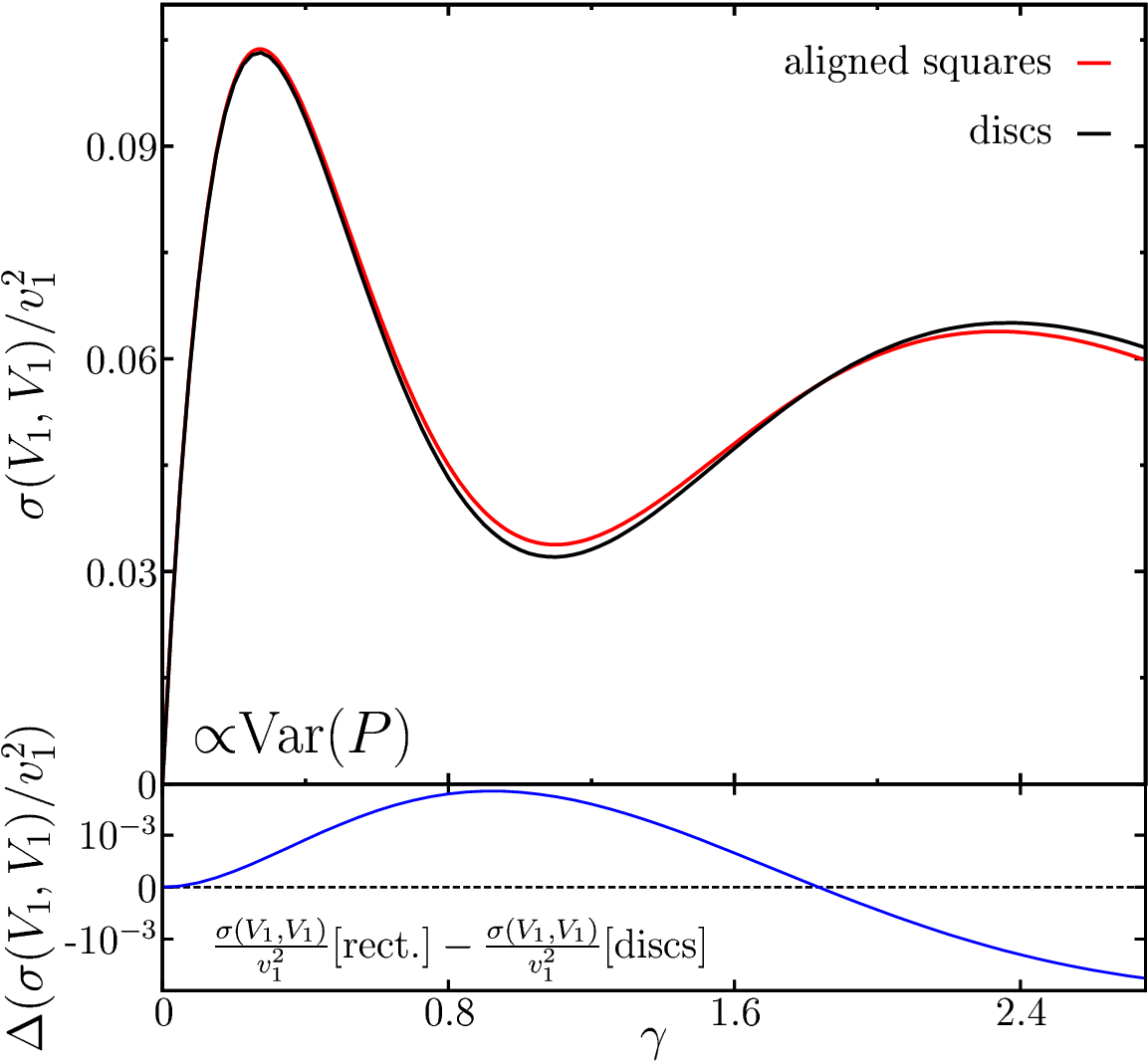}
  \hfill
  \includegraphics[width=0.45\textwidth]{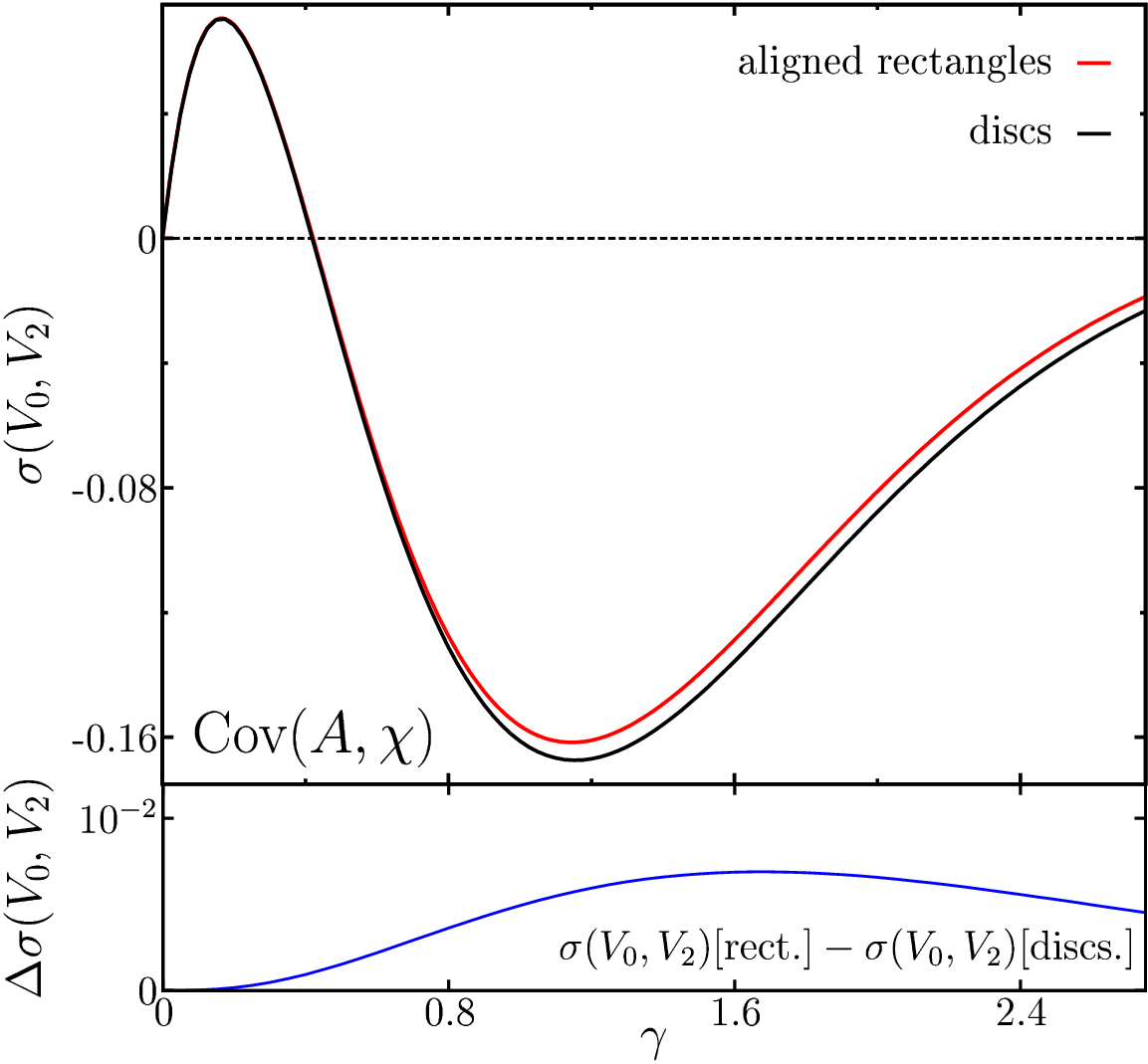}
  \hfill\hbox{}\\ \hfill
  \includegraphics[width=0.45\textwidth]{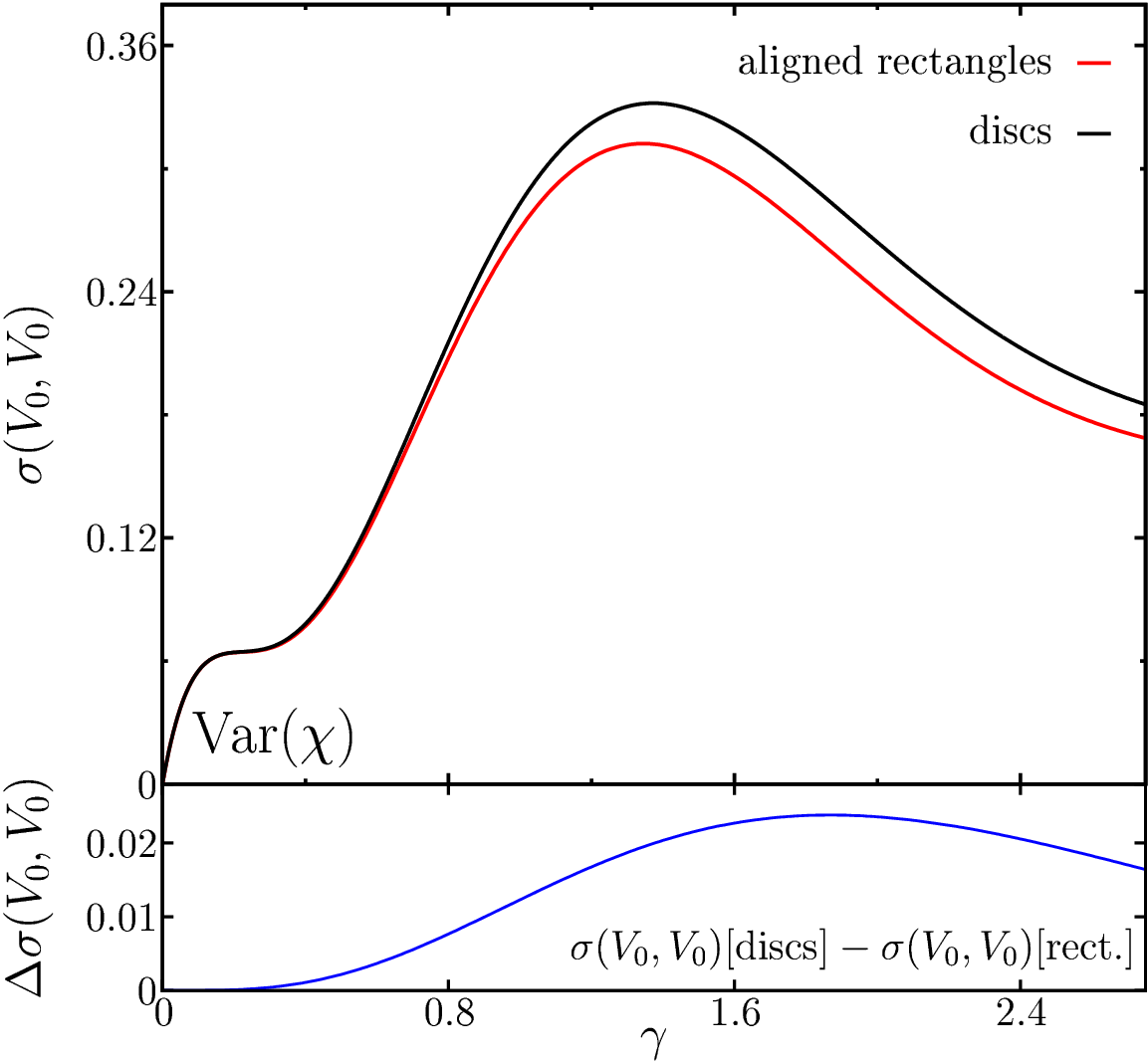}
  \hfill
  \includegraphics[width=0.45\textwidth]{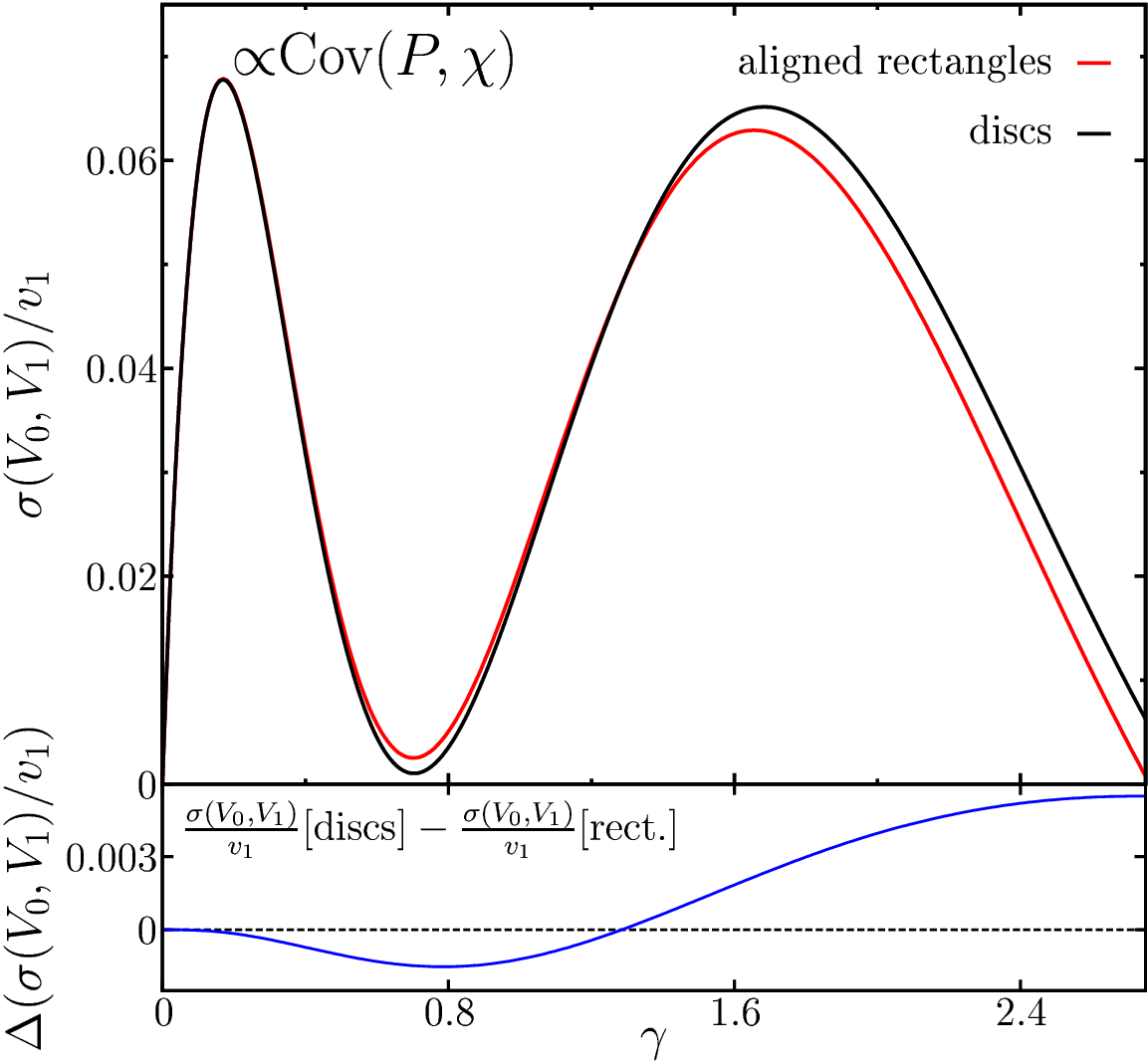}
  \hfill\hbox{} 
  \caption{Variances and covariances: $\sigma(V_2,V_2)$,
    the variance of the area; $\sigma(V_1,V_2)/v_1$, proportional to the
    covariance of area and perimeter; $v_1$ is half of the perimeter
    of a single grain; $\sigma(V_1,V_1)/v_1^2$, proportional to the
    variance of the perimeter; $\sigma(V_0,V_2)$, the covariance of
    area and Euler characteristic; $\sigma(V_0,V_0)$, the variance of
    the Euler characteristic; $\sigma(V_0,V_1)/v_1$, proportional to the
    covariance of perimeter and Euler characteristic. They are shown
    both for Boolean models with aligned rectangles, see Eqs.~\protect\eqref{rect2},
    \protect\eqref{rect5}, \protect\eqref{4522}, \protect\eqref{4527}, \protect\eqref{sigma01}, and
    \protect\eqref{rect6789}, and for overlapping discs, see~\cite{HLS}. Note
    that except for $\sigma(V_1,V_1)/v_1^2$
    the curves for the rectangles are independent of the aspect ratio of
    the rectangle, see also Figs.~\protect\ref{fig_cov_numerical_integration}
    and \protect\ref{fig_Boolean_cov}. The insets in the figures at the top are
    close-up views which show that the covariances differ slightly for Boolean
    model with rectangles or with discs. Below each subfigure, the differences
    of the covariances for Boolean models with rectangles or discs are plotted.}
    \label{fig_bm_cov_rect_discs}
\end{figure}

The variances and covariances of the Minkowski functionals of overlapping
rectangles exhibit a complex behaviour as functions of the intensity
$\gamma$ similar to the Boolean model with discs in~\cite[Section 7]{HLS}. The
variances of area and Euler characteristic apparently have one maximum and no other
extrema. The variance of the perimeter has a global maximum and (at
least) one local minimum. As
expected, the three Minkowski functionals are positively correlated at
low intensities $\gamma$, but at relatively high intensities the area is
anti-correlated to both the Euler characteristic and the perimeter.

The covariance $\sigma(V_0,V_1)$ between the perimeter and the Euler
characteristic shows a qualitatively different behaviour for rectangles
with an isotropic orientation distribution when compared to the overlapping
discs or aligned rectangles.
There is a regime in the
intensity $\gamma$ (around the first local minimum) for which the
rectangles with an isotropic orientation distribution are
anti-correlated, while the aligned grains are positively correlated like
the discs in~\cite{HLS}. This is probably related to the fact that
rotated rectangles can more easily form clusters with holes than aligned
rectangles or discs. The zero-crossing of the expectation of the Euler
characteristic $\chi$ for the rectangles with an isotropic orientation
distribution is within this regime. For the aligned rectangles, the
zero-crossing of the mean value of $\chi$ is at the end of this
regime, see~\cite{HoerrmannHugKlattMecke}.

The question remains whether or not the variances and covariances of
area or rescaled perimeter of the Boolean model are independent of the
grain distribution like the first moments of these functionals.
Equations~\eqref{rect2} and \eqref{rect5} show that at least for aligned
rectangles the variance $\sigma(V_2,V_2)$ as well as the covariance
$\sigma(V_1,V_2)$ divided by the perimeter of a single grain $(2a+2b)$
are indeed independent of the aspect ratio.
The simulation results from Fig.~\ref{fig_Boolean_cov} might suggest
that this could also be valid for the isotropic orientation
distributions.
However, the variance $\sigma(V_2,V_2)$ and the rescaled covariance
$\sigma(V_1,V_2)/(2a+2b)$ do depend on the grain distribution, although
only weakly for the models studied here. To show this, we evaluate
Eqs.~\eqref{rect2} and \eqref{rect5} numerically and compare the
covariances to those of the Boolean model with discs from~\cite{HLS}.
Figure~\ref{fig_bm_cov_rect_discs} shows that there is a weak but
significant difference in the analytic curves of $\sigma(V_2,V_2)$ and
$\sigma(V_1,V_2)$ for the two different
models.
The variance of the perimeter depends more clearly on the grain
distribution. Even if it is rescaled by the perimeter of a single grain
and even for aligned grains, the variance distinctly depends on the
aspect ratio of the rectangles (except for small intensities $\gamma$).
So, in contrast to the first moments of the area and rescaled perimeter
of the Boolean model, the second moments in general depend on the grain
distribution, e.g., the orientation distribution, even if this
dependence may be weak. 
As expected, also the variance $\sigma(V_0,V_0)$ of the
Euler characteristic as well as the covariances $\sigma(V_0,V_1)$ and
$\sigma(V_0,V_2)$ depend on the grain distribution, see
Fig.~\ref{fig_bm_cov_rect_discs}.

\subsection{Central limit theorem}
\label{sec_Boolean_epdf}

\begin{figure}[t] \centering
  \subfigure{\includegraphics[width=0.3\textwidth]{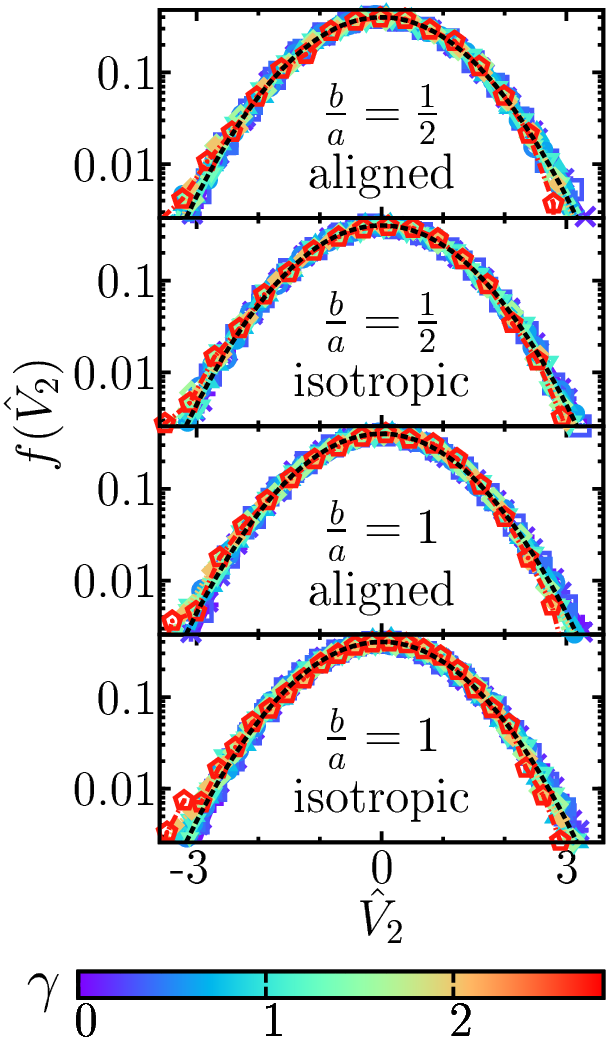}}%
  \hspace*{2pt}
  \subfigure{\includegraphics[width=0.3\textwidth]{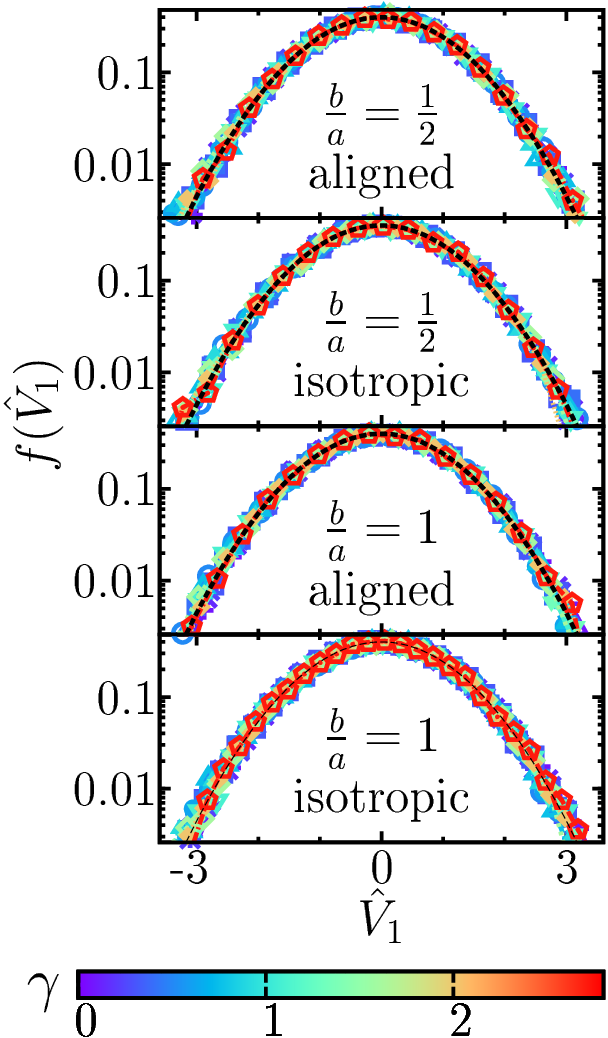}}%
  \hspace*{2pt}
  \subfigure{\includegraphics[width=0.3\textwidth]{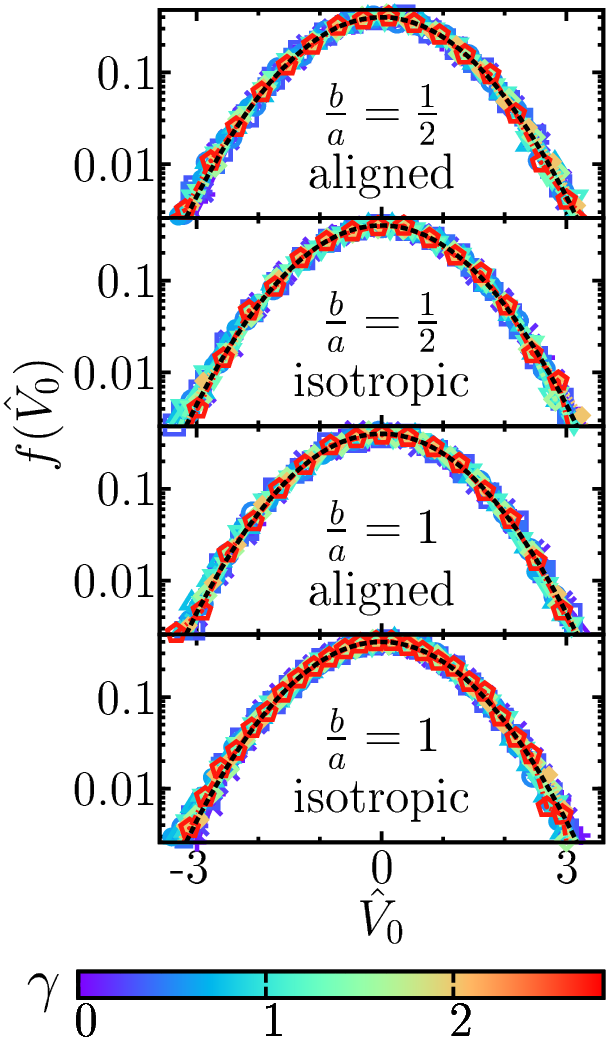}}%
  \caption{Histograms $f$ of
    the normalized intrinsic volumes $\hat{V}_i$ (see
    Eq.~\protect\eqref{eq:normalizedVi}) of Boolean models of rectangles with
    different aspect ratios $b/a$ for different intensities $\gamma$;
    for both aligned rectangles and rectangles with an isotropic
    orientation distribution. For all of these different models, the
    rescaled distributions are already for the relatively small system
  size $L=20a$ in very good agreement with a normal distribution (dashed
black line).}
\label{fig_Boolean_MF_epdf}
\end{figure}

We also determine the histograms of the intrinsic volumes in a finite
observation window, where the histograms are weighted by the total
number of samples and the bin width. The histograms are then compared to
the density of a standard normal distribution in order to numerically validate
the central limit theorems in Section~\ref{sec:CLTs}.
The information content of a histogram is up to the binning almost
equivalent to the empirical distribution function, but in plotting it is
more convenient to compare histograms and densities.

The histograms resemble probability density functions. However,
the intrinsic volumes of the considered Boolean model do not have
probability density functions. Indeed, with positive probability, there
is no overlap between the grains and there are no intersections with the
boundary so that some multiples of the intrinsic volumes of the fixed
grain have positive probability.

In this subsection, we simulate larger systems than in the previous
Subsection~\ref{sec:var-cov-simulation}. For a simulation box
with side length $L=20a$, we perform for each different set of
parameters between $M_s=\num{5000}$ and \num{150000} simulations of
Boolean models with rectangles at varying intensities $\gamma$.
Like in Subsection~\ref{sec:var-cov-simulation}, the rectangles have aspect
ratio $1$ or $1/2$, and they are either aligned w.r.t. the observation
window or their orientation is isotropically distributed.
To produce the histograms, we simulate more than \num{1400000} samples
of Boolean models including about \num{350000000} rectangles in total.

We normalize the intrinsic volumes $V_i$, i.e., we subtract the
estimated mean values $\langle{}V_i{}\rangle$ of the intrinsic volumes
and divide by $\sqrt{s(V_i,V_i)}$:
\begin{align}
  \hat{V}_i := \frac{V_i - \langle{}V_i{}\rangle}{\sqrt{s(V_i,V_i)}}\;.
  \label{eq:normalizedVi}
\end{align}
Figure~\ref{fig_Boolean_MF_epdf} plots the histograms $f$ of the
normalized intrinsic volumes of Boolean models with different aspect
ratios $b/a$ for either aligned rectangles or rectangles with an
isotropic orientation distribution and for varying intensities $\gamma$.

{These histograms are in good agreement with the density function of a normal
distribution for all intrinsic volumes, for all intensities, and for all
of the simulated models (despite the relatively small simulation box).}
In other words, even in small observation windows, the probability
distributions of the intrinsic volumes of Boolean models can be well
approximated by Gaussian distributions.

As we have mentioned in the introduction, the central limit theorems for
the geometric functionals (see Theorems~\ref{thm:MultivariateCase} and
\ref{thm:UnivariateCase}) and the exact formulae for the second moments
(see Theorems~\ref{thm:CovariancesVolumeSurfaceArea} and
\ref{thm:PositiveDefiniteness}) can be used for hypothesis testing of
models of random heterogeneous media.
A hypothesis test could, e.g., use the intrinsic volumes to decide whether
or not a random two-phase medium can be modeled by overlapping grains.
The joint probability distribution of the Minkowski functionals allows
for a characterization of the shape by several geometrical
functionals and hence for a construction of tests using their full
covariance structure. For a different random field (with a Poisson distributed
number of counts in a binned gamma-ray sky map) such a sensitive
morphometric data analysis has already been
developed~\cite{GoeringKlattEtAl2013, KlattGoeringStegmannMecke2012}.
The same concepts could be applied to the Boolean model.

\begin{figure}[t]
  \centering
  \includegraphics[width=0.9\textwidth]{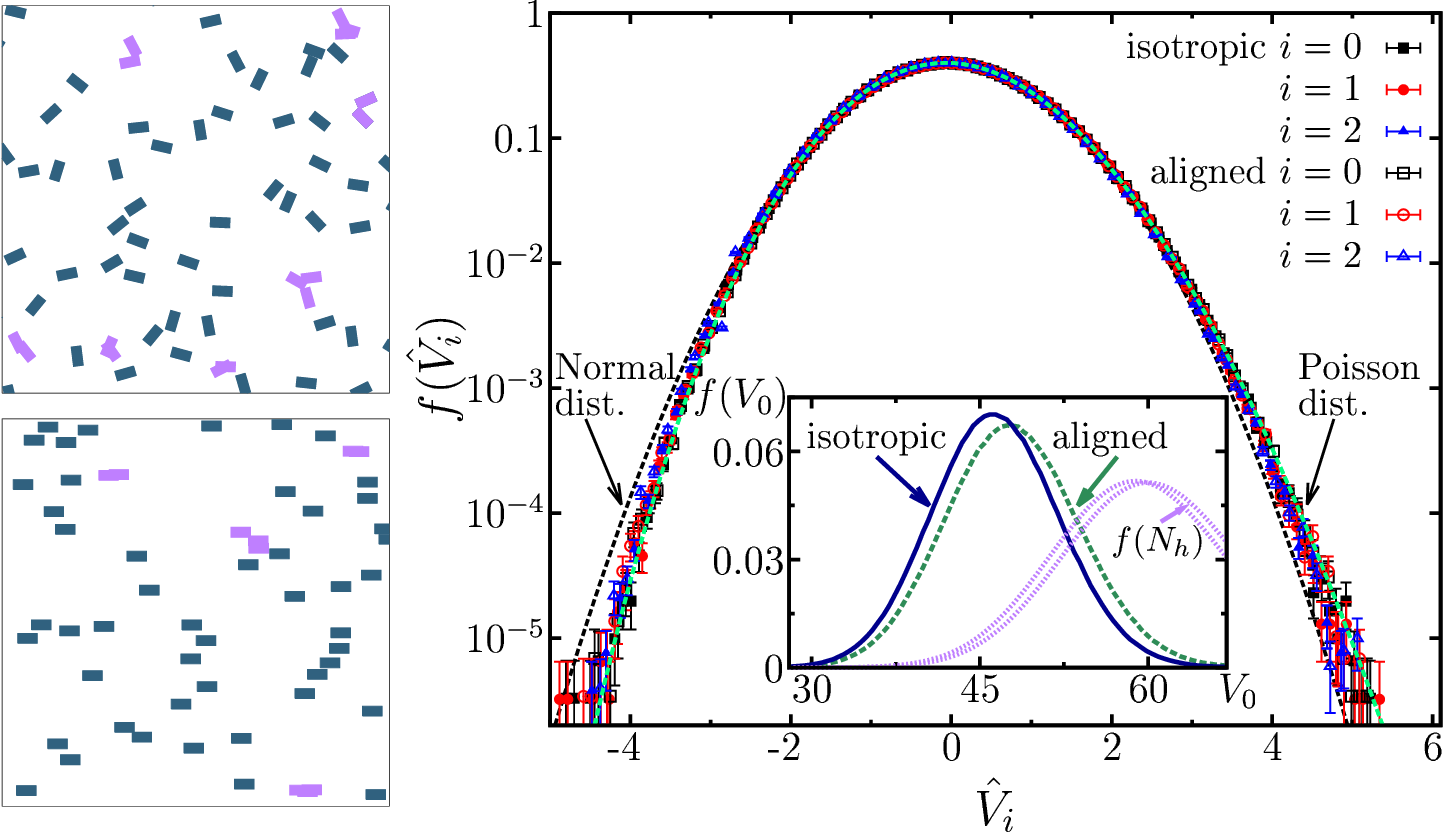}
  \caption{Histograms $f$ of the normalized intrinsic volumes $\hat{V}_i$ (see
    Eq.~\protect\eqref{eq:normalizedVi})
    of Boolean models with either aligned
    rectangles or rectangles with an isotropic orientation distribution.
    The dashed black line depicts the density of a normal distribution.
    On the left
    hand side, two samples of the Boolean models with
    either an isotropic orientation distribution or aligned rectangles
    are shown; clusters of rectangles are colored purple.
    The inset shows the histograms of the non-rescaled Euler characteristic
    $V_0$ in the isotropic and aligned case as well as the
    corresponding probability mass functions of the number
    of grains $N_h$ hitting the observation window with mean values
    $59.4$ and $60.5$, respectively.}
    \label{fig_Boolean_MF_epdf_detail}
\end{figure}

{In Fig.~\ref{fig_Boolean_MF_epdf}, there are only small deviations from
a normal distribution relative to the error bars. So, the systematic
deviations, e.g., due to the finite observation window size, seem to be
small.}
In order to determine these deviations, a very high numerical accuracy
is needed.
We simulate $3\cdot 10^6$ samples of two Boolean models for rectangles
with aspect ratio ${1}/{2}$ that are either aligned or follow an
isotropic orientation distribution. Here, we apply minus sampling
boundary conditions, i.e., we consider all grains with centres in a
slightly larger simulation box
$[-\sqrt{a^2+b^2}/2,L+\sqrt{a^2+b^2}/2]^2$, but the observation window
is still the original square $(0,L)^2$ with ${L}=20a$. 
Contributions caused by the boundary are here neglected as it is often
done in physics.
The expected number of grains in the simulation box is adjusted
accordingly and follows a Poisson distribution with mean $\gamma\cdot
(L+\sqrt{a^2+b^2})^2$.
To minimize the computational costs, a relatively low intensity $\gamma$
is chosen for these simulations. It corresponds to an expected occupied
area fraction $\phi={1}/{15}$.
The resulting histograms are plotted in
Fig.~\ref{fig_Boolean_MF_epdf_detail}. 
As expected, the high statistics reveal for the small system size
deviations from the normal distribution that are significantly larger
than the error bars.

For each underlying Boolean model, the histograms of the normalized
intrinsic volumes coincide within error bars. This is not surprising
because at the low intensity chosen here the intrinsic volumes are
strongly correlated. (The correlation coefficients are larger than 0.9.) 

For different Boolean models (isotropic orientation distribution or
aligned grains), the histogram of the non-rescaled Minkowski functionals differ
slightly but distinctly already for the relatively small intensity
studied here, see the inset of Fig.~\ref{fig_Boolean_MF_epdf_detail}.
In contrast to this, the histograms of the normalized intrinsic
volumes collapse for the different Boolean models within the error bars
to a single curve, which can be well approximated by a standardized
Poisson distribution.
This can be explained by the strong correlation between the intrinsic
volumes and the number of grains $N_h$ hitting the observation window
for each Boolean model. The latter follows a Poisson distribution with
parameter $\mathbf{E}[N_h] = \gamma \cdot \mathbf{E}[V_2([0;L]^2+Z_0)]$. 
(The correlation coefficients are larger than 0.85.)
There is only a small relative difference between the parameters
$\mathbf{E}[N_h]$ for the different considered Boolean models, because
the observation window is large when compared to the typical grain $Z_0$.
Therefore, the corresponding Poisson distributions are very close after
standardization (dashed green line in
Fig.~\ref{fig_Boolean_MF_epdf_detail}) and coincide with the histograms
of the normalized intrinsic volumes. 

\begin{acknowledgement}
 We would like to thank Julia H\"orrmann and Klaus Mecke for some
 helpful remarks and discussions. The authors were supported by the
 German Research Foundation (DFG) through the research unit ``Geometry
 and Physics of Spatial Random Systems''.
\end{acknowledgement}

\end{document}